\documentclass{amsart}
\usepackage{amsmath,amssymb}
\usepackage[dvips]{graphics}
\newtheorem{Theorem}{Theorem}[section]
\newtheorem{Lemma}[Theorem]{Lemma}
\newtheorem{Proposition}[Theorem]{Proposition}

\theoremstyle{Definition}

\newtheorem{Example}[Theorem]{Example}

\theoremstyle{Remark}
\newtheorem{Remark}[Theorem]{Remark}

\def\labelenumi{\rm ({\makebox[0.8em][c]{\theenumi}})}
\def\theenumi{\arabic{enumi}}
\def\labelenumii{\rm {\makebox[0.8em][c]{(\theenumii)}}}
\def\theenumii{\roman{enumii}}

\renewcommand{\thefigure}
{\arabic{section}.\arabic{figure}}
\makeatletter
\@addtoreset{figure}{section}%{subsection}
\def\@thmcountersep{-}
\makeatother

\numberwithin{equation}{section}

\begin{document} 

\title{Homotopy on spatial graphs and the Sato-Levine invariant}

\author{Thomas Fleming}
\address{Department of Mathematics, University of California San Diego, 9500 Gilman Drive, La Jolla, CA 92093, USA}
%\curraddr{}
\email{tfleming@math.ucsd.edu}
\thanks{The first author was supported by a Fellowship of 
the Japan Society for the Promotion of Science for Post-Doctoral Foreign 
Researchers (Short-Term) (No. PE05003).}

\author{Ryo Nikkuni}
\address{Institute of Human and Social Sciences, Faculty of Teacher Education, Kanazawa University, Kakuma-machi, Kanazawa, Ishikawa, 920-1192, Japan}
\email{nick@ed.kanazawa-u.ac.jp}
\thanks{The second author was partially supported by a Grant-in-Aid for Scientific Research (B) (2) (No. 15340019), Japan Society for the Promotion of Science.}

%    General info

\date{}

\dedicatory{}


\begin{abstract}
Edge-homotopy and vertex-homotopy are equivalence relations on spatial 
graphs which are generalizations of Milnor's link-homotopy. 
We introduce some edge (resp. vertex)-homotopy invariants of spatial 
graphs by applying the Sato-Levine invariant for the $2$-component 
constituent algebraically split links and show examples 
of non-splittable spatial graphs up to edge (resp. vertex)-homotopy, 
all of whose constituent links are link-homotopically trivial. 
\end{abstract}

\maketitle

\section{Introduction} 

Throughout this paper we work in the piecewise linear category. 
Let $G$ be a finite graph which does not have isolated vertices and 
free vertices. An embedding $f$ of $G$ into the $3$-sphere $S^{3}$ is called  
a {\it spatial embedding of $G$} or simply a {\it spatial graph}. 
For a spatial embedding $f$ and a subgraph $H$ of $G$ which is homeomorphic 
to the $1$-sphere $S^{1}$ or a disjoint union of $1$-spheres, we call $f(H)$ 
a {\it constituent knot} or a 
{\it constituent link} of $f$, respectively. 
A graph $G$ is said to be {\it planar} if there exists an embedding of $G$ 
into the $2$-sphere $S^{2}$, and a spatial embedding of a planar 
graph is said to be {\it trivial} if it is ambient isotopic to an embedding 
of the graph into a $2$-sphere in $S^{3}$. 
A spatial embedding $f$ of a graph $G$ is said to be 
{\it split} if there exists a $2$-sphere $S$ in $S^{3}$ such that 
$S\cap f(G)=\emptyset$ and each component of $S^{3}-S$ 
has intersection with $f(G)$, and otherwise $f$ is said to be 
{\it non-splittable}. 

Two spatial embeddings of a graph $G$ are said to be {\it edge-homotopic} if 
they are transformed into each other by {\it self crossing changes} and 
ambient isotopies, where a self crossing change is a crossing change on 
the same spatial edge, and {\it vertex-homotopic} if they are transformed 
into each other by crossing changes on two adjacent spatial edges and 
ambient isotopies.\footnote{In \cite{taniyama94b}, edge-homotopy and 
vertex-homotopy were called {\it homotopy} and 
{\it weak homotopy}, respectively.} 
These equivalence relations were introduced by Taniyama 
\cite{taniyama94b} as generalizations of Milnor's 
{\it link-homotopy} on links \cite{milnor54}, namely if $G$ is 
homeomorphic to a disjoint union of $1$-spheres, then these are 
none other than link-homotopy. 
There are many studies about link-homotopy. In particular, 
the link-homotopy classification was given for $2$- and 
$3$-component links by Milnor \cite{milnor54}, for $4$-component links by 
Levine \cite{levine88} and for all links by 
Habegger and Lin \cite{habegger-lin90}. On the other hand, 
there are very few studies about edge (resp. vertex)-homotopy on 
spatial graphs \cite{taniyama94a}, 
\cite{motohashi-taniyama96}, \cite{ohyama-taniyama01}, \cite{nikkuni04}.  

In \cite{taniyama94a}, Taniyama defined an edge 
(resp. vertex)-homotopy invariant of spatial graphs called the 
{\it $\alpha$-invariant} by applying the {\it Casson invariant} 
(or equivalently 
the second coefficient of the {\it Conway polynomial}) of the constituent 
knots and showed that there exists a non-trivial spatial embedding $f$ of 
a planar graph up to edge (resp. vertex)-homotopy, 
even in the case where $f$ does 
not contain any constituent link. 
But the $\alpha$-invariant cannot detect a non-splittable spatial 
embedding of a disconnected graph up to edge (resp. vertex)-homotopy. 
As far as the authors know, an example 
of a non-splittable spatial embedding of 
a disconnected graph up to edge (resp. vertex)-homotopy, all of whose 
constituent links are link-homotopically trivial has not yet been 
demonstrated.  

Our purpose in this paper is to study spatial embeddings 
of disconnected graphs up to edge (resp. vertex)-homotopy by applying 
the {\it Sato-Levine invariant} \cite{sato84} (or equivalently the third 
coefficient of the Conway polynomial) for the constituent 
$2$-component algebraically split links and show that there 
exist infinitely many non-splittable spatial embeddings of a certain 
disconnected graph up to edge (resp. vertex)-homotopy all of whose 
constituent links are link-homotopically trivial.  
These examples show that edge (resp. vertex)-homotopy on spatial graphs 
behaves quite differently than link-homotopy on links. 
In the next 
section we give the definitions of our invariants and state their 
invariance up to edge (resp. vertex)-homotopy. 

\section{Definitions of invariants} 

We call a subgraph of a graph $G$ a {\it cycle} if it 
is homeomorphic to the $1$-sphere, and a cycle is called a 
{\it $k$-cycle} if it contains exactly $k$ edges. 
For a subgraph $H$ of $G$, 
we denote the set of all cycles of $H$ by $\Gamma(H)$. 
We set ${\mathbb Z}_{m}=\{0,1,\ldots,m-1\}$ 
for a positive 
integer $m$ and ${\mathbb Z}_{0}={\mathbb Z}$.  
We regard ${\mathbb Z}_{m}$ as an abelian group in the obvious way. 
We call a map $\omega:\Gamma(H)\to {\mathbb Z}_{m}$ a 
{\it weight on $\Gamma(H)$ over ${\mathbb Z}_{m}$}. For an edge $e$ of $H$, 
we denote the set of all cycles of $H$ which contain the edge $e$ by 
$\Gamma_{e}(H)$. For a pair of two adjacent edges 
$e_{1}$ and $e_{2}$ of $H$, we denote the set of all 
cycles of $H$ which contain 
the edges $e_{1}$ and $e_{2}$ by $\Gamma_{e_{1},e_{2}}(H)$. 
Then we say that a weight $\omega$ on $\Gamma(H)$ over ${\mathbb Z}_{m}$ is 
{\it weakly balanced \footnote{
A weight $\omega$ on $\Gamma(H)$ over ${\mathbb Z}_{m}$ 
is said to be {\it balanced on an edge $e$ of $H$} if 
$\sum_{\gamma\in\Gamma_{e}(H)}\omega(\gamma)[\gamma]=0$ in 
$H_{1}(H;{\mathbb Z}_{m})$, where the orientation of 
$\gamma$ is induced by the one of $e$ \cite{taniyama94a}.} 
on an edge $e$} if 
\begin{eqnarray*}
\sum_{\gamma\in\Gamma_{e}(H)}\omega(\gamma)=0
\end{eqnarray*}
in ${\mathbb Z}_{m}$ \cite{nikkuni02}, and {\it weakly balanced on a pair of 
adjacent edges $e_{1}$ and $e_{2}$} if 
\begin{eqnarray*}
\sum_{\gamma\in\Gamma_{e_{1},e_{2}}(H)}\omega(\gamma)=0
\end{eqnarray*}
in ${\mathbb Z}_{m}$. 
Let $G=G_{1}\cup G_{2}$ be a disjoint union of two connected graphs and 
$\omega_{i}:\Gamma(G_{i})\to {\mathbb Z}_{m}$ a weight on $\Gamma(G_{i})$ 
over ${\mathbb Z}_{m}$ $(i=1,2)$. Let $f$ be a spatial embedding of $G$ 
such that 
\begin{eqnarray*}
\omega_{1}(\gamma)\omega_{2}(\gamma')
{\rm lk}(f(\gamma),f(\gamma'))=0
\end{eqnarray*}
in ${\mathbb Z}$ for any $\gamma\in \Gamma(G_{1})$ and 
$\gamma'\in \Gamma(G_{2})$, where 
${\rm lk}(L)={\rm lk}(K_{1},K_{2})$ denotes the {\it linking number} 
of a $2$-component oriented link $L=K_{1}\cup K_{2}$. 
Then we define $\beta_{\omega_{1},\omega_{2}}(f)\in {\mathbb Z}_{m}$ by 
\begin{eqnarray*}
\beta_{\omega_{1},\omega_{2}}(f)\equiv \sum_{\gamma\in\Gamma(G_{1}) \atop \gamma'\in\Gamma(G_{2})}
\omega_{1}(\gamma)\omega_{2}(\gamma')
a_{3}(f(\gamma),f(\gamma'))\pmod{m}, 
\end{eqnarray*}
where $a_{3}(L)=a_{3}(K_{1},K_{2})$ denotes the third coefficient of the 
Conway polynomial of a $2$-component oriented link $L=K_{1}\cup K_{2}$. 
We remark here that $a_{3}(L)$ coincides with 
the {\it Sato-Levine invariant} $\beta(L)$ of 
$L$ if $L$ is {\it algebraically split}, namely ${\rm lk}(K_{1},K_{2})=0$ 
\cite{cochran85}, \cite{sturm90}. 
Thus our $\beta_{\omega_{1},\omega_{2}}(f)$ 
is also the modulo $m$ reduction of the summation 
of Sato-Levine invariants for the constituent $2$-component algebraically 
split links of $f$. 

\begin{Remark}\label{sl_rem}
{\rm For a $2$-component algebraically split link $L=K_{1}\cup K_{2}$,  
\begin{enumerate}
\item The value of $a_{3}(L)$ does not depend on the orientations 
of $K_{1}$ and $K_{2}$. 
Actually we can check it easily by the original definition of the 
Sato-Levine invariant. 
\item The value of $a_{3}(L)$ is not a link-homotopy invariant of $L$ 
(see also Lemma \ref{keylemma}). For example, 
the Whitehead link $L$ is link-homotopically trivial but $a_{3}(L)=1$. 
\end{enumerate}
}
\end{Remark}

Now we state the invariance of $\beta_{\omega_{1},\omega_{2}}$ up to 
edge (resp. vertex)-homotopy under some conditions on the graphs. 

\begin{Theorem}\label{inv1}
Let $G=G_{1}\cup G_{2}$ be a disjoint union of two connected graphs and 
$\omega_{i}$ a weight on $\Gamma(G_{i})$ 
over ${\mathbb Z}_{m}$ $(i=1,2)$. 
Let $f$ be a spatial embedding of $G$ such that 
\begin{eqnarray*}
\omega_{1}(\gamma)\omega_{2}(\gamma')
{\rm lk}(f(\gamma),f(\gamma'))=0
\end{eqnarray*}
in ${\mathbb Z}$ for any $\gamma\in \Gamma(G_{1})$ and $\gamma'\in \Gamma(G_{2})$. 
Then we have the following: 
\begin{enumerate}
\item If $\omega_{i}$ is weakly balanced on any edge of $G_{i}$ $(i=1,2)$, 
then $\beta_{\omega_{1},\omega_{2}}(f)$ is an edge-homotopy invariant of $f$. 
\item If $\omega_{i}$ is weakly balanced on any pair of adjacent edges 
of $G_{i}$ $(i=1,2)$, then 
$\beta_{\omega_{1},\omega_{2}}(f)$ is a vertex-homotopy invariant of $f$. 
\end{enumerate}
\end{Theorem}

We prove Theorem \ref{inv1} in the next section. In addition, by using an 
integer-valued invariant (Theorem \ref{integer_inv}), we show that there 
exist infinitely many non-splittable spatial embeddings of a certain 
disconnected graph up to edge-homotopy all of whose constituent links are 
link-homotopically trivial (Example \ref{integer_inv_ex}). 
We also exhibit an infinite family of non-splittable spatial embeddings of 
a certain disconnected graph up to vertex-homotopy which can be 
distinguished by our integer-valued invariant (Example \ref{integer_inv_ex2}). 

\vspace{0.2cm}
We note that if a graph $G$ contains a connected component which is 
homeomorphic to the $1$-sphere, then our invariants in Theorem \ref{inv1} 
are useless. For such cases, we can define edge (vertex)-homotopy 
invariants that take values in ${\mathbb Z}_{2}$ on weaker condition for 
weights than the one stated in Theorem \ref{inv1}. 
For a subgraph $H$ of a graph $G$, 
we say that a weight $\omega$ on $\Gamma(H)$ over ${\mathbb Z}_{2}$ is 
{\it totally balanced} if 
\begin{eqnarray*}
\sum_{\gamma\in\Gamma(H)}\omega(\gamma)[\gamma]=0
\end{eqnarray*}
in $H_{1}(H;{\mathbb Z}_{2})$. We note that if a weight 
$\omega$ on $\Gamma(H)$ over ${\mathbb Z}_{2}$ is totally balanced, then 
it is weakly balanced on any edge $e$ of $H$ (Lemma \ref{keysublemma}), 
but not always weakly 
balanced on any pair of adjacent edges of $H$ 
(Remark \ref{keysublemmaremark}). 
Then we have the following. 

\begin{Theorem}\label{mod2inv}
Let $G=G_{1}\cup G_{2}$ be a disjoint union of two connected graphs 
and $\omega_{i}$ a weight on $\Gamma(G_{i})$ 
over ${\mathbb Z}_{2}$ $(i=1,2)$. 
Let $f$ be a spatial embedding of $G$ such that 
\begin{eqnarray*}
\omega_{1}(\gamma)\omega_{2}(\gamma'){\rm lk}(f(\gamma),f(\gamma'))=0
\end{eqnarray*}
in ${\mathbb Z}$ for any $\gamma\in \Gamma(G_{1})$ and $\gamma'\in \Gamma(G_{2})$. 
Then we have the following: 
\begin{enumerate}
\item If either $\omega_{1}$ is totally balanced 
on $\Gamma(G_{1})$ 
or $\omega_{2}$ is totally balanced 
on $\Gamma(G_{2})$, then 
$\beta_{\omega_{1},\omega_{2}}(f)$ is an edge-homotopy invariant of $f$. 
\item If either $\omega_{1}$ is totally balanced on $\Gamma(G_{1})$ and 
weakly balanced on any pair of adjacent edges 
of $G_{1}$, or 
$\omega_{2}$ is totally balanced on 
$\Gamma(G_{2})$ and weakly balanced on any pair of adjacent edges 
of $G_{2}$, then 
$\beta_{\omega_{1},\omega_{2}}(f)$ is a vertex-homotopy invariant of $f$. 
\end{enumerate}
\end{Theorem}

We also prove Theorem \ref{mod2inv} in the next section and give some 
examples in Section $5$. In particular, we show that there 
exist infinitely many non-splittable spatial embeddings of a certain 
disconnected graph up to vertex-homotopy, all of whose constituent 
links are link-homotopically trivial (Example \ref{mod2inv_ex2}). 
We remark here that the ${\mathbb Z}_{2}$-valued invariant in 
Theorem \ref{mod2inv} cannot always be extended to an integer-valued 
one (Remark \ref{mod2inv_rem}). 

\begin{figure}[htbp]
      \begin{center}
\scalebox{0.4}{\includegraphics*{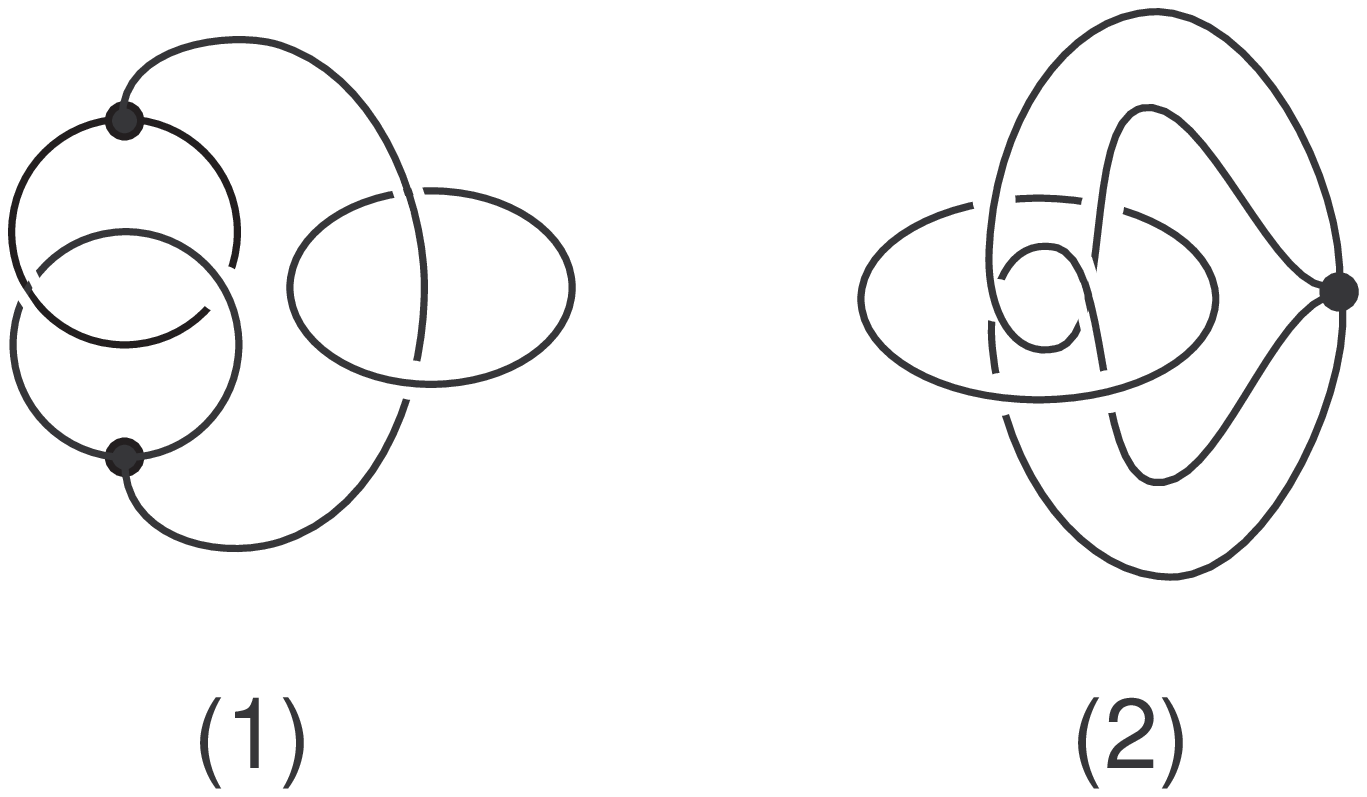}}
      \end{center}
   \caption{}
  \label{handcuff_bouquet}
\end{figure}

\vspace{0.2cm}
Theorems \ref{inv1} and \ref{mod2inv} do not work for spatial graphs 
as illustrated in Figure \ref{handcuff_bouquet}, for instance. In 
Section $6$, we state a method to detect such non-splittable spatial 
graphs up to edge-homotopy by using a planar surface having a graph 
as a spine (Theorem \ref{modulo2_inv_more}). 
Actually we show that each of the spatial graphs 
as illustrated in Figure \ref{handcuff_bouquet} is non-splittable up to 
edge-homotopy (Example \ref{hand_bouquet_ex}). 

%
%\begin{figure}[htbp]
%      \begin{center}
%  \epsfile{file=handcuff_bouquet,scale=0.4} 
%      \end{center}
%   \caption{}
%  \label{handcuff_bouquet}
%\end{figure} 
% 
%
\section{Proofs of Theorems \ref{inv1} and \ref{mod2inv}} 

We first calculate the change in the third coefficient of the Conway 
polynomial of $2$-component algebraically 
split links which differ by a single self crossing change. 

\begin{figure}[htbp]
      \begin{center}
\scalebox{0.4}{\includegraphics*{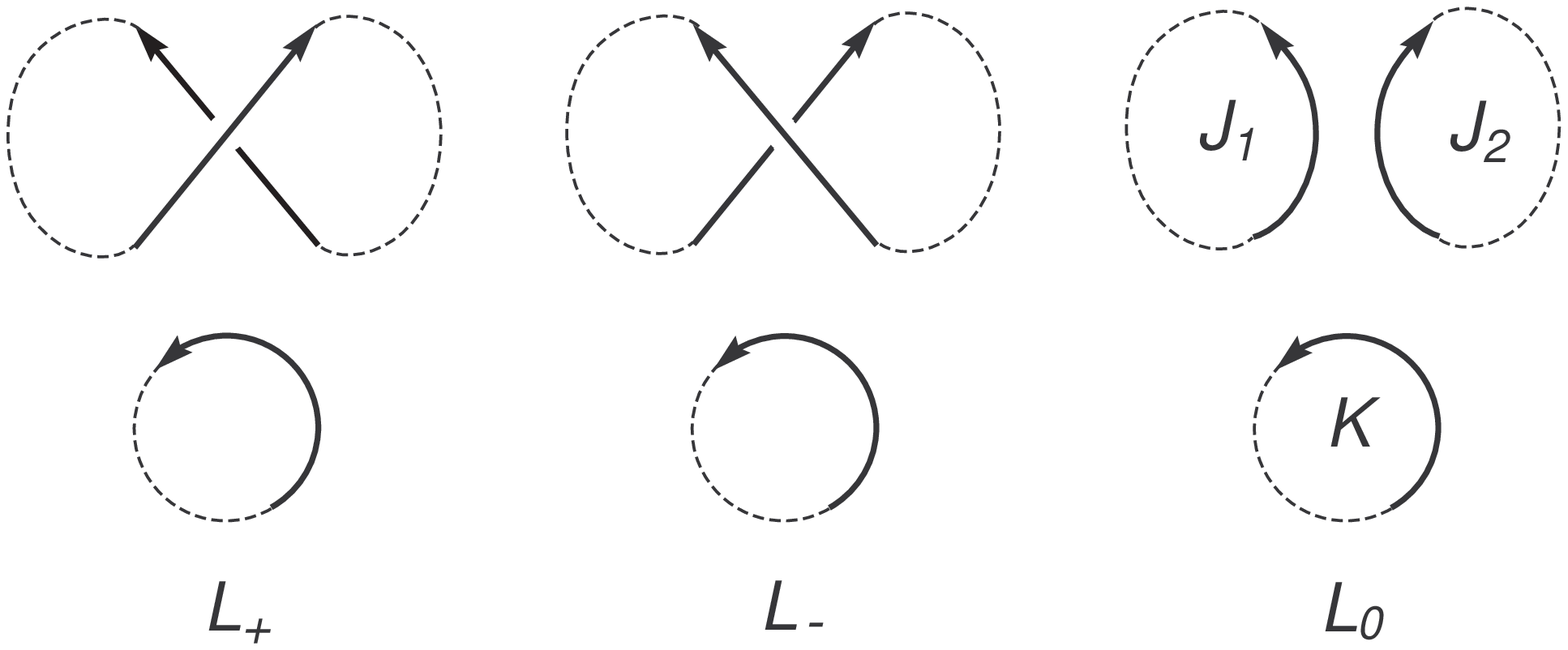}}
      \end{center}
   \caption{}
  \label{skein_sc}
\end{figure}

\begin{Lemma}\label{keylemma}
Let $L_{+}$ and $L_{-}$ be two $2$-component oriented links and 
$L_{0}=J_{1}\cup J_{2}\cup K$ a $3$-component oriented link which 
are identical except inside the depicted regions as illustrated in Figure 
\ref{skein_sc}. Suppose that ${\rm lk}(L_{+})={\rm lk}(L_{-})=0$. 
Then it holds that 
\begin{eqnarray*}
a_{3}(L_{+})-a_{3}(L_{-})=-{\rm lk}(J_{1},K)^{2}=-{\rm lk}(J_{2},K)^{2}. 
\end{eqnarray*}
\end{Lemma}

\begin{proof} 
By the skein relation of the Conway polynomial and 
a well-known formula for 
the second coefficient of the Conway 
polynomial of a $3$-component oriented 
link (cf. \cite{hosokawa58}, \cite{hartley83}, \cite{hoste85}), we have that 
\begin{eqnarray}
a_{3}(L_{+})-a_{3}(L_{-})
&=&{\rm lk}(J_{1},J_{2}){\rm lk}(J_{2},K)+{\rm lk}(J_{2},K){\rm lk}(J_{1},K)
\label{skein1} \\
&&+{\rm lk}(J_{1},K){\rm lk}(J_{1},J_{2}). \nonumber
\end{eqnarray}
We note that 
\begin{eqnarray}\label{lk}
{\rm lk}(J_{1},K)+{\rm lk}(J_{2},K)=0
\end{eqnarray}
by the condition ${\rm lk}(L_{+})={\rm lk}(L_{-})=0$. Thus by 
(\ref{skein1}) and (\ref{lk}), we have that 
\begin{eqnarray*}
a_{3}(L_{+})-a_{3}(L_{-})&=&
{\rm lk}(J_{1},J_{2})\left\{-{\rm lk}(J_{1},K)\right\}+{\rm lk}(J_{2},K){\rm lk}(J_{1},K)\\
&&+{\rm lk}(J_{1},K){\rm lk}(J_{1},J_{2})\\
&=&{\rm lk}(J_{2},K){\rm lk}(J_{1},K). 
\end{eqnarray*}
Therefore by (\ref{lk}) we have the result. 
\end{proof}

\begin{figure}[htbp]
      \begin{center}
\scalebox{0.4}{\includegraphics*{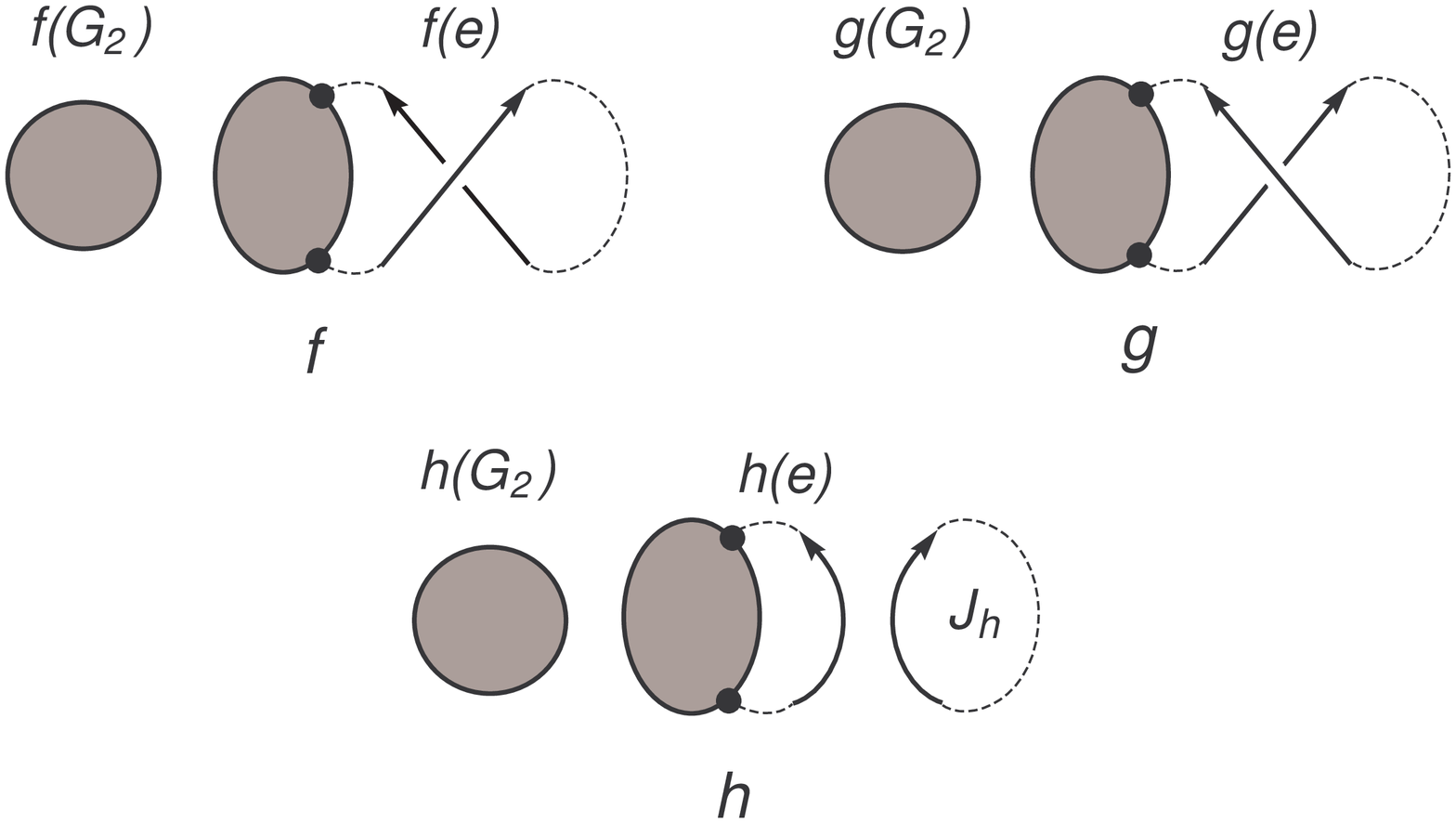}}
      \end{center}
   \caption{}
  \label{skein_sc2}
\end{figure}

\begin{figure}[htbp]
      \begin{center}
\scalebox{0.4}{\includegraphics*{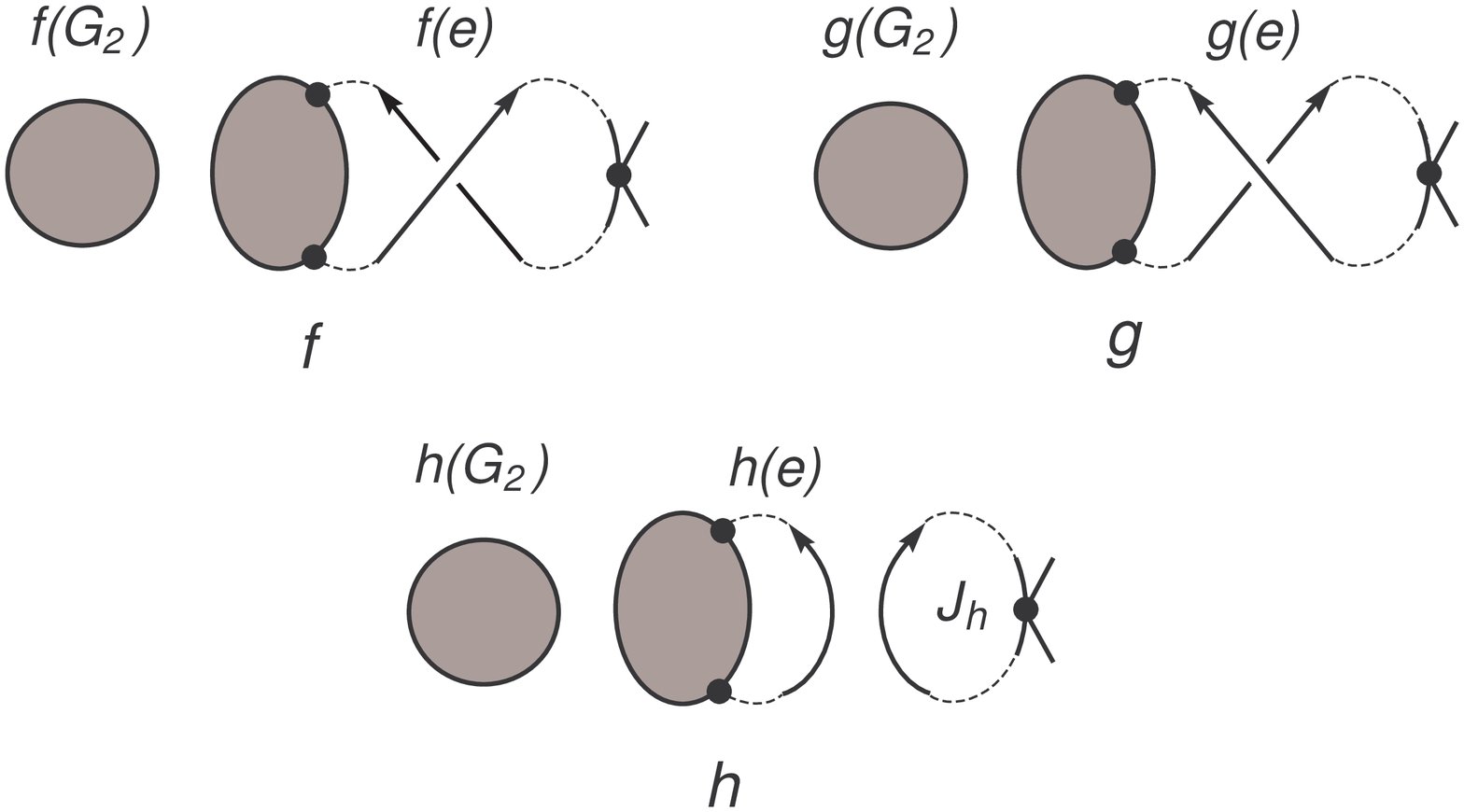}}
      \end{center}
   \caption{}
  \label{skein_sc3}
\end{figure}

\begin{proof}[Proof of Theorem \ref{inv1}.]
(1) Let $f$ and $g$ be two spatial embeddings of $G$ such that 
\begin{eqnarray}\label{lk0f}
\omega_{1}(\gamma)\omega_{2}(\gamma')
{\rm lk}(f(\gamma),f(\gamma'))
=0
\end{eqnarray}
in ${\mathbb Z}$ for any $\gamma\in \Gamma(G_{1})$ and $\gamma'\in \Gamma(G_{2})$ 
and $g$ is edge-homotopic to $f$. 
Then it also holds that 
\begin{eqnarray}\label{lk0g}
\omega_{1}(\gamma)\omega_{2}(\gamma')
{\rm lk}(g(\gamma),g(\gamma'))=0
\end{eqnarray}
in ${\mathbb Z}$ for any $\gamma\in \Gamma(G_{1})$ and 
$\gamma'\in \Gamma(G_{2})$ because the linking number of a 
$2$-component constituent link of a 
spatial graph is an edge-homotopy invariant. 
First we show that if $f$ is transformed into $g$ by 
self crossing changes on $f(G_{1})$ and ambient isotopies, then 
$\beta_{\omega_{1},\omega_{2}}(f)=\beta_{\omega_{1},\omega_{2}}(g)$. 
It is clear that any link invariant of a constituent link of a spatial graph 
is also an ambient isotopy invariant of the spatial graph. Thus 
we may assume that $g$ is obtained from $f$ by a single crossing change 
on $f(e)$ for an edge $e$ of $G_{1}$ 
as illustrated in Figure \ref{skein_sc2}. Moreover, by smoothing this 
crossing point we can obtain the spatial embedding $h$ 
of $G$ and the knot $J_{h}$ as illustrated in Figure \ref{skein_sc2}. 
Then by (\ref{lk0f}), (\ref{lk0g}), Lemma \ref{keylemma} and 
the assumption for $\omega_{1}$ 
we have that 
\begin{eqnarray*}
\beta_{\omega_{1},\omega_{2}}(f)-\beta_{\omega_{1},\omega_{2}}(g)
&\equiv& \sum_{\gamma\in\Gamma(G_{1}) \atop \gamma'\in\Gamma(G_{2})}
\omega_{1}(\gamma)\omega_{2}(\gamma')
\left\{a_{3}(f(\gamma),f(\gamma'))-a_{3}(g(\gamma),g(\gamma'))\right\}\\
&=&\sum_{\gamma\in\Gamma_{e}(G_{1}) \atop \gamma'\in\Gamma(G_{2})}
\omega_{1}(\gamma)\omega_{2}(\gamma')
\left\{a_{3}(f(\gamma),f(\gamma'))-a_{3}(g(\gamma),g(\gamma'))\right\}\\
&=&-\sum_{\gamma\in\Gamma_{e}(G_{1}) \atop \gamma'\in\Gamma(G_{2})}
\omega_{1}(\gamma)\omega_{2}(\gamma')
{\rm lk}(h(\gamma'),J_{h})^{2}\\
&=&-\left(\sum_{\gamma\in \Gamma_{e}(G_{1})}\omega_{1}(\gamma)\right)
\sum_{\gamma'\in\Gamma(G_{2})}\omega_{2}(\gamma')
{\rm lk}(h(\gamma'),J_{h})^{2}\\
&\equiv&0.  
\end{eqnarray*}
Therefore we 
have that 
$\beta_{\omega_{1},\omega_{2}}(f)=\beta_{\omega_{1},\omega_{2}}(g)$.  
In the same way we can show that if $f$ is transformed into $g$ by 
self crossing changes on $f(G_{2})$ and ambient isotopies, then 
$\beta_{\omega_{1},\omega_{2}}(f)=\beta_{\omega_{1},\omega_{2}}(g)$. 
Thus we have that $\beta_{\omega_{1},\omega_{2}}$ is an 
edge-homotopy invariant. 

(2) By considering 
the triple of spatial embeddings as illustrated in Figure \ref{skein_sc3}, 
we can prove (2) in a similar way as the proof of (1). 
We omit the details. 
\end{proof}


Next we prove Theorem \ref{mod2inv}. For a subgraph $H$ of a graph $G$, 
we have the following. 

\begin{Lemma}\label{keysublemma}
A totally balanced weight $\omega$ on $\Gamma(H)$ over ${\mathbb Z}_{2}$ is 
weakly 
balanced on any edge $e$ of $H$. 
\end{Lemma}

\begin{proof} 
For an edge $e$ of $H$, we can represent any 
$\gamma\in \Gamma_{e}(H)$ as $e+c_{\gamma}\in Z_{1}(H;{\mathbb Z}_{2})$, where 
$c_{\gamma}$ is a $1$-chain in $C_{1}(H\setminus{e};{\mathbb Z}_{2})$. 
Then we have that 
\begin{eqnarray*}
0&=&\sum_{\gamma\in\Gamma(H)}\omega(\gamma)[\gamma]\\
&=&\sum_{\gamma\in\Gamma_{e}(H)}\omega(\gamma)[e+c_{\gamma}]
+\sum_{\gamma'\in\Gamma(H)\setminus\Gamma_{e}(H)}\omega(\gamma')[\gamma']
\end{eqnarray*}
in $H_{1}(H;{\mathbb Z}_{2})$. 
This implies that if $\omega$ is not weakly 
balanced on $e$, then $\omega$ is not totally balanced 
on $\Gamma(H)$ over ${\mathbb Z}_{2}$. 
\end{proof}
\begin{Remark}\label{keysublemmaremark}
{\rm 
A totally balanced weight $\omega$ on $\Gamma(H)$ over ${\mathbb Z}_{2}$ is 
not always weakly 
balanced on any pair of adjacent edges of $H$. 
For example, 
let $\omega$ be a weight on $\Theta_{3}$ 
(see Example \ref{integer_inv_ex}) 
over ${\mathbb Z}_{2}$ 
defined by $\omega(\gamma)=1$ for any cycle 
$\gamma\in \Gamma(\Theta_{3})$. It is easy to see 
that $\omega$ is totally balanced, but not weakly 
balanced, on each pair of adjacent edges of $\Theta_{3}$.}
\end{Remark}

\begin{proof}[Proof of Theorem \ref{mod2inv}.] 
(1) Let $f$ and $g$ be two spatial embeddings of $G$ which are edge-homotopic 
such that 
\begin{eqnarray*}\label{lk02}
\omega_{1}(\gamma)\omega_{2}(\gamma')
{\rm lk}(f(\gamma),f(\gamma'))
=
\omega_{1}(\gamma)\omega_{2}(\gamma')
{\rm lk}(g(\gamma),g(\gamma'))
=0
\end{eqnarray*}
in ${\mathbb Z}$ for any $\gamma\in \Gamma(G_{1})$ and $\gamma'\in \Gamma(G_{2})$. 
First we show that if $f$ is transformed into $g$ by 
self crossing changes on $f(G_{1})$ and ambient isotopies, then 
$\beta_{\omega_{1},\omega_{2}}(f)=\beta_{\omega_{1},\omega_{2}}(g)$. 
In the same way as the proof of Theorem \ref{inv1}, we may consider 
three spatial embeddings $f,g$ and $h$ of $G$ and the knot $J_{h}$ as 
illustrated in Figure \ref{skein_sc2}. 
Then, by the same calculation in the proof of Theorem \ref{inv1}, 
we have that 
\begin{eqnarray*}
\beta_{\omega_{1},\omega_{2}}(f)-\beta_{\omega_{1},\omega_{2}}(g)
&\equiv& 
-\left(\sum_{\gamma\in \Gamma_{e}(G_{1})}\omega_{1}(\gamma)\right)
\sum_{\gamma'\in\Gamma(G_{2})}\omega_{2}(\gamma')
{\rm lk}(h(\gamma'),J_{h})^{2}\\
&\equiv&  \left(\sum_{\gamma\in \Gamma_{e}(G_{1})}\omega_{1}(\gamma)\right)
\sum_{\gamma'\in\Gamma(G_{2})}\omega_{2}(\gamma')
{\rm lk}(h(\gamma'),J_{h})\\
&\equiv&\left(\sum_{\gamma\in \Gamma_{e}(G_{1})}\omega_{1}(\gamma)\right)
{\rm lk}\left(\sum_{\gamma'\in\Gamma(G_{2})}
\omega_{2}(\gamma')h(\gamma'),J_{h}\right). 
\end{eqnarray*}
If $\omega_{1}$ is totally balanced on $\Gamma(G_{1})$, then by 
Lemma \ref{keysublemma} it is 
weakly balanced on any edge $e$ of $G_{1}$. This implies that 
$\beta_{\omega_{1},\omega_{2}}(f)=\beta_{\omega_{1},\omega_{2}}(g)$. 
If $\omega_{2}$ is totally 
balanced on $\Gamma(G_{1})$, then we have that 
\begin{eqnarray*}
{\rm lk}\left(\sum_{\gamma'\in\Gamma(G_{2})}
\omega_{2}(\gamma')h(\gamma'),J_{h}\right)\equiv
{\rm lk}\left(0,J_{h}\right)=0.
\end{eqnarray*}
Therefore this also implies that $\beta_{\omega_{1},\omega_{2}}(f)
=\beta_{\omega_{1},\omega_{2}}(g)$. 
In the same way we can show that if $f$ is transformed into $g$ by 
self crossing changes on $f(G_{2})$ and ambient isotopies, then 
$\beta_{\omega_{1},\omega_{2}}(f)=\beta_{\omega_{1},\omega_{2}}(g)$. 
Thus we have that $\beta_{\omega_{1},\omega_{2}}$ is an 
edge-homotopy invariant. 

(2) By considering 
the triple of spatial embeddings as illustrated in Figure \ref{skein_sc3}, 
we can prove (2) in a similar way as the proof of (1). 
We also omit the details. 
\end{proof}

\vspace{0.3cm}
Since the Conway polynomial of a split link is zero, 
our invariants take the value zero for any split ($2$-component) 
spatial graph. Therefore if the value of our invariant of a spatial graph 
is not zero, then it is non-splittable up to edge (resp. vertex)-homotopy.

\section{Integer-valued invariants} 

\begin{figure}[htbp]
      \begin{center} 
\scalebox{0.4}{\includegraphics*{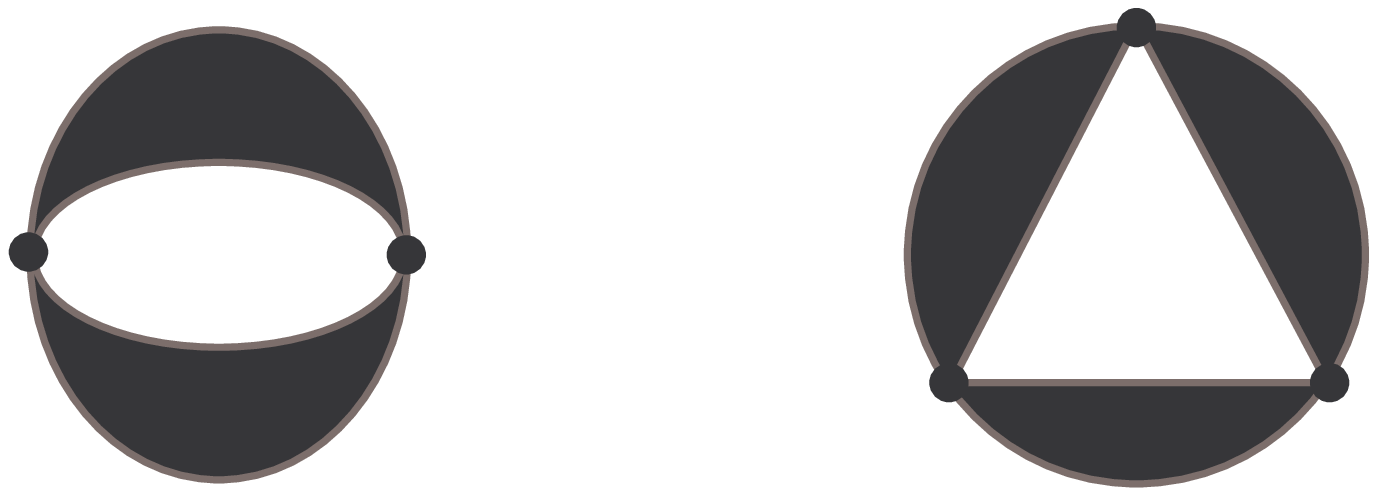}}
      \end{center}
   \caption{}
  \label{checkerboard}
\end{figure}

Let $G$ be a planar graph. An embedding $p:G\to S^{2}$ is said to be 
{\it cellular} if the closure of each of the connected components of 
$S^{2}-p(G)$ is homeomorphic to the disk. Then we regard the set of 
the boundaries of all of the connected 
components of $S^{2}-p(G)$ as a subset of $\Gamma(G)$ and denote it 
by $\Gamma_{p}(G)$. 
We say that {\it $G$ admits a 
checkerboard coloring on $S^{2}$} if there exists a cellular embedding  
$p:G\to S^{2}$ such that we can color all of the connected components 
of $S^{2}-p(G)$ by two colors (black and white) so that any of the two 
components which are adjacent by an edge have distinct colors; see 
Figure \ref{checkerboard}. 
We denote the subset of $\Gamma_{p}(G)$ which corresponds to the black 
(resp. white) colored components by $\Gamma_{p}^{b}(G)$ 
(resp. $\Gamma_{p}^{w}(G)$). 

\begin{Proposition}\label{integer_inv_weight}
Let $G$ be a planar graph which is not homeomorphic to $S^{1}$ and 
admits a checkerboard coloring on $S^{2}$ 
with respect to a cellular embedding $p:G\to S^{2}$. 
Let $\omega_{p}$ be a weight on $\Gamma(G)$ over 
${\mathbb Z}$ defined by 
\begin{eqnarray*}
\omega_{p}(\gamma)= \left\{
\begin{array}{@{\,}ll}
1 & \mbox{$(\gamma\in \Gamma_{p}^{b}(G))$,} \\
-1 &   \mbox{$(\gamma\in\Gamma_{p}^{w}(G))$,} \\
0 &   \mbox{$(\gamma\in \Gamma(G)\setminus \Gamma_{p}(G))$.}
\end{array}
\right. 
\end{eqnarray*}
Then $\omega_{p}$ is weakly balanced on any edge of $G$. 
\end{Proposition}

\begin{proof} 
For any edge $e$ of $G$, there exist exactly two cycles 
$\gamma\in \Gamma_{p}^{b}(G)$ and $\gamma'\in \Gamma_{p}^{w}(G)$ such that 
$e\subset\gamma$ and $e\subset\gamma'$. Thus we have the result. 
\end{proof}

We call the weight $\omega_{p}$ in Proposition \ref{integer_inv_weight} a 
{\it checkerboard weight}. Thus by 
Proposition \ref{integer_inv_weight} and Theorem \ref{inv1} (1), 
we can obtain 
an integer-valued edge-homotopy invariant as follows. 

\begin{Theorem}\label{integer_inv}
Let $G=G_{1}\cup G_{2}$ be a disjoint union of two connected planar 
graphs such that $G_{i}$ is not homeomorphic to $S^{1}$ and 
admits a checkerboard coloring on $S^{2}$ 
with respect to a cellular embedding $p_{i}:G_{i}\to S^{2}$ $(i=1,2)$. 
Let $\omega_{p_{i}}$ be a checkerboard 
weight on $\Gamma(G_{i})$ over ${\mathbb Z}$ $(i=1,2)$ 
and $f$ a spatial embedding of $G$ such that 
\begin{eqnarray*}
\omega_{p_{1}}(\gamma)\omega_{p_{2}}(\gamma')
{\rm lk}(f(\gamma),f(\gamma'))=0
\end{eqnarray*}
in ${\mathbb Z}$ for any $\gamma\in \Gamma(G_{1})$ and 
$\gamma'\in \Gamma(G_{2})$. 
Then $\beta_{\omega_{p_{1}},\omega_{p_{2}}}(f)$ is an integer-valued 
edge-homotopy invariant of $f$. \hfill $\square$
\end{Theorem}
\begin{Example}\label{integer_inv_ex}
{\rm 
Let $\Theta_{n}$ be a graph with two vertices $u$ and $v$ and $n$ edges 
$e_{1},e_{2},\ldots,e_{n}$, each of which joins $u$ and $v$. A spatial 
embedding of $\Theta_{n}$ is called a {\it (spatial) theta $n$-curve} 
or simply a {\it theta curve} if $n=3$. For $n\ge 2$, 
we denote that a cycle of $\Theta_{n}$ 
consists of two edges $e_{i}$ and $e_{j}$ by $\gamma_{ij}$ ($i<j$). 
Then it is clear that $\Theta_{n}$ admits a cellular embedding 
$p:\Theta_{n}\to S^{2}$ so that 
\begin{eqnarray*}
\Gamma_{p}(\Theta_{n})
=\left\{\gamma_{12},\gamma_{23},\ldots,\gamma_{n-1,n},\gamma_{1n}\right\}.
\end{eqnarray*}
Moreover, for $m\ge 1$, $\Theta_{2m}$ admits a checkerboard coloring on 
$S^{2}$ so that 
\begin{eqnarray*}
\Gamma_{p}^{b}(\Theta_{2m})
&=&\left\{\gamma_{12},\gamma_{34},\ldots,\gamma_{2m-1,2m}\right\}, \\
\Gamma_{p}^{w}(\Theta_{2m})
&=&\left\{\gamma_{23},\gamma_{45},\ldots,\gamma_{2m-2,2m-1},\gamma_{1,2m}\right\}. 
\end{eqnarray*}

Now let $G$ be a disjoint union of two copies of $\Theta_{4}$, 
each of which admits a checkerboard coloring on $S^{2}$ with respect to 
the cellular embedding $p$ as above. 
Let $\omega_{p}$ be a checkerboard weight on $\Gamma(\Theta_{4})$ 
over ${\mathbb Z}$ and $g_{1}$ a spatial embedding of $G$ as illustrated 
in Figure \ref{theta4_theta4}. We can see that any of the 
$2$-component constituent 
links of $g_{1}$ has a zero linking number. 
More precisely, 
$g_{1}$ contains exactly one 
non-trivial $2$-component link 
$L=g_{1}(\gamma_{14})\cup g_{1}(\gamma_{14}')$ whose linking number is zero. 
Thus by Theorem \ref{integer_inv} we have that 
$\beta_{\omega_{p},\omega_{p}}(g_{1})$ is an integer-valued edge-homotopy 
invariant of $g_{1}$. 
Then, by a direct calculation we have that $a_{3}(L)=2$,  
namely $\beta_{\omega_{p},\omega_{p}}(g_{1})=2$. 
Note that a $2$-component link is link-homotopically trivial if and only if 
its linking number is zero \cite{milnor54}. 
This implies 
that $g_{1}$ is non-splittable up to edge-homotopy despite the fact that 
any of the constituent links of $g_{1}$ is link-homotopically trivial.  

\begin{figure}[htbp]
      \begin{center}
\scalebox{0.35}{\includegraphics*{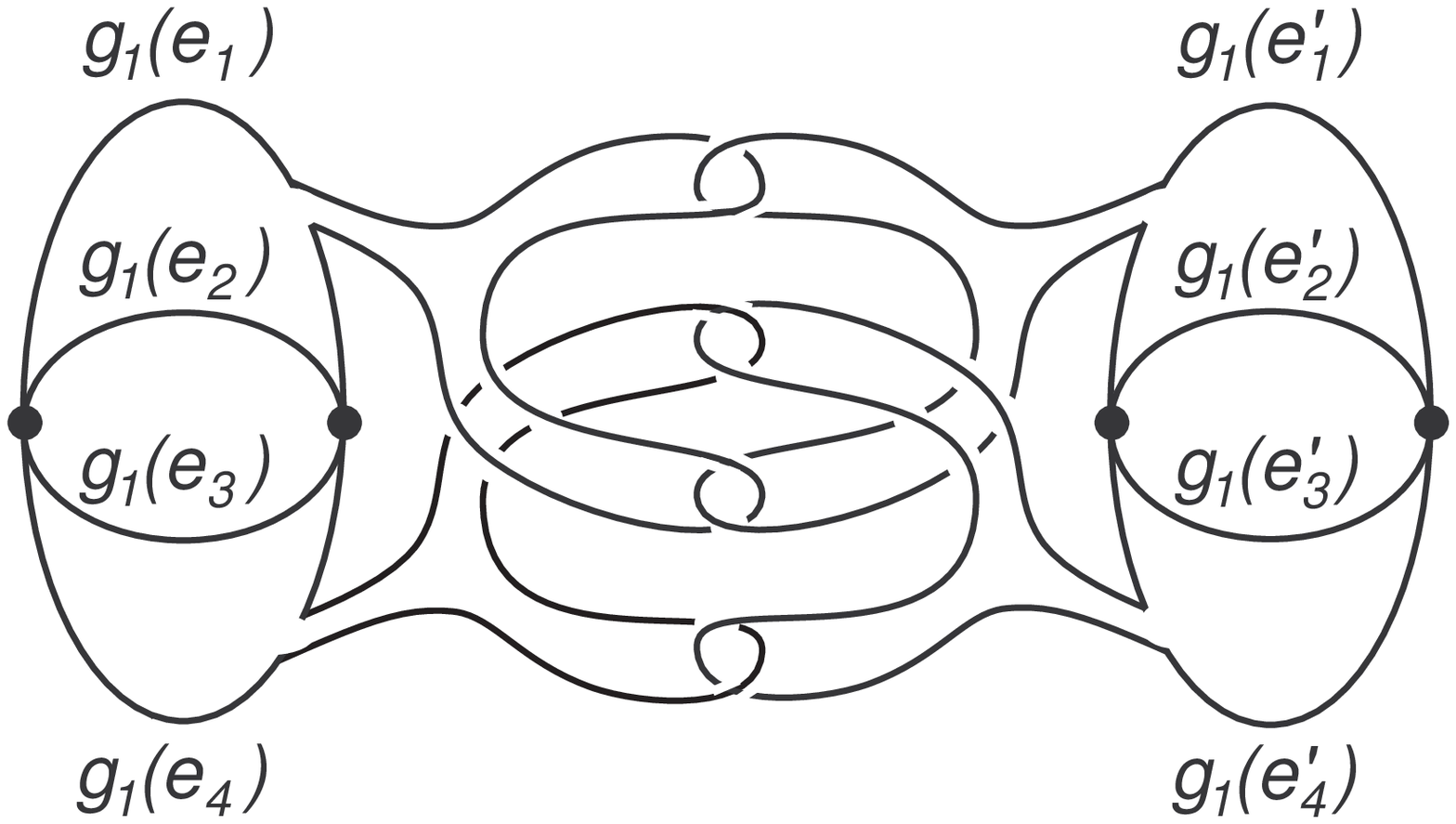}}
      \end{center}
   \caption{}
  \label{theta4_theta4}
\end{figure} 

Moreover, for an integer $m$, let $g_{m}$ be a spatial embedding of $G$ 
as illustrated in Figure \ref{bandsum_C3}. 
If $m\neq 0$, we can see that $g_{m}$ contains exactly one 
non-trivial $2$-component link 
$L=g_{m}(\gamma_{14})\cup g_{m}(\gamma_{14}')$ whose linking number is zero. 
Thus we also have that 
$\beta_{\omega_{p},\omega_{p}}(g_{m})$ is an integer-valued edge-homotopy 
invariant of $g_{m}$. 
Then, by a calculation we have that $a_{3}(L)=2m$,  
namely $\beta_{\omega_{p},\omega_{p}}(g_{m})=2m$. This implies 
that there exist infinitely many non-splittable spatial embeddings of $G$ 
up to edge-homotopy, all of whose constituent links are link-homotopically 
trivial. 

\begin{figure}[htbp]
      \begin{center}
\scalebox{0.35}{\includegraphics*{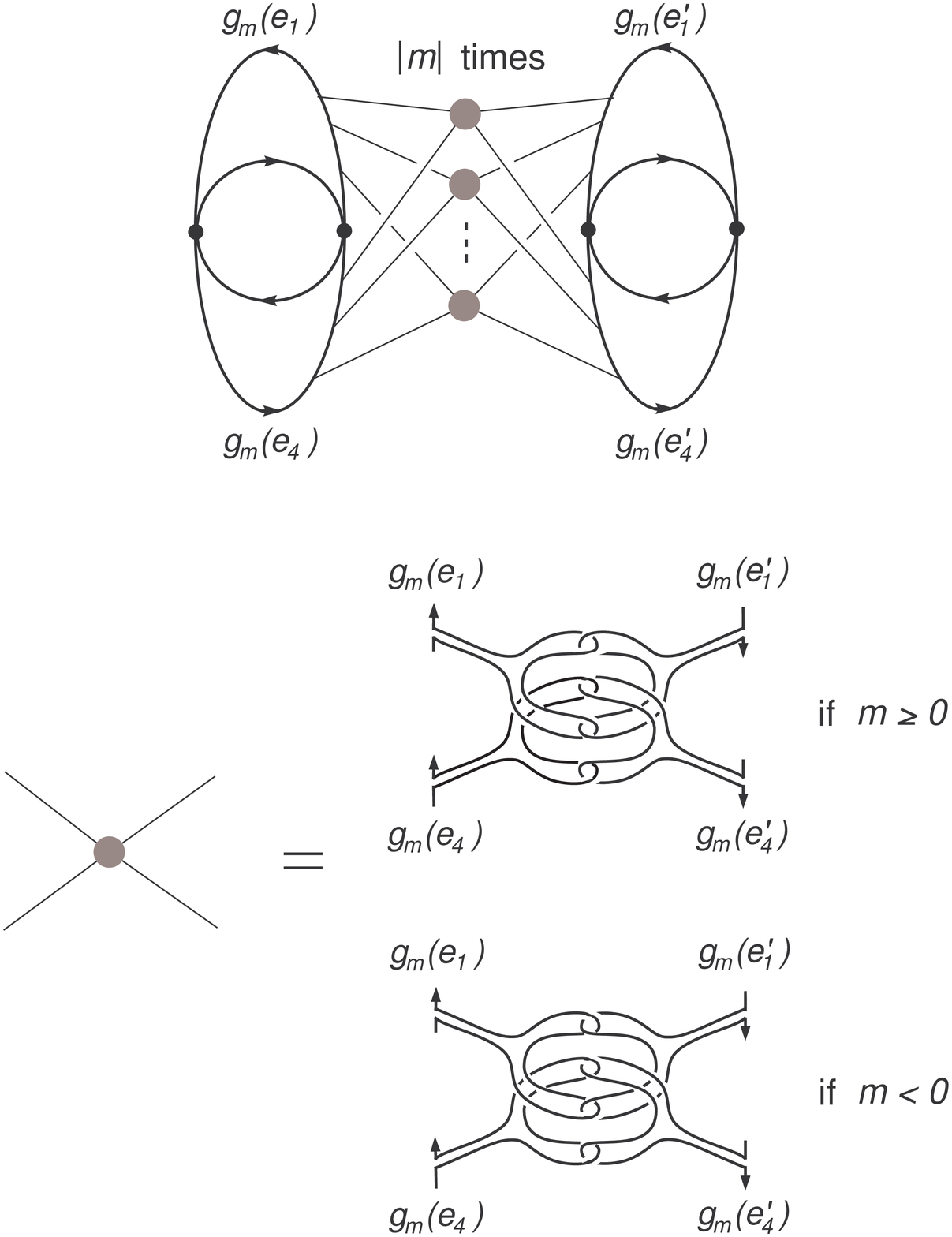}}
      \end{center}
   \caption{}
  \label{bandsum_C3}
\end{figure} 
}
\end{Example}
\begin{Example}\label{integer_inv_ex2}
{\rm 
Let $H$ be a graph as illustrated in Figure \ref{VH_ex_graph}. 
We denote the cycle of $H$ which contains $e_{i}$ and $e_{j}$ 
by $\gamma_{ij}$ ($i<j$). Let $G$ be a disjoint union of two copies of 
$H$ and $g_{1}$ a spatial embedding of $G$ as illustrated in 
Figure \ref{VH_ex}. This spatial embedding $g_{1}$ contains exactly one 
$4$-component constituent link 
$L=g_{1}(\gamma_{12}\cup \gamma_{34}\cup \gamma'_{12}\cup \gamma'_{34})$. 
Note that if $g_{1}$ is split up to vertex-homotopy, then $L$ is split 
up to link-homotopy. Since $|\mu_{1234}(L)|=1$, where $\mu_{1234}$ denotes 
Milnor's {\it $\mu$-invariant} of length $4$ of $4$-component links 
\cite{milnor54}, we have that $L$ is non-splittable up to link-homotopy. 
Therefore we have that $g_{1}$ is non-splittable up to vertex-homotopy. 

We can also prove this fact by our integer-valued vertex-homotopy invariant 
as follows. Let $\omega$ be a weight 
on $\Gamma(H)$ over ${\mathbb Z}$ 
defined by $\omega(\gamma_{14})=\omega(\gamma_{23})=1$, 
$\omega(\gamma_{13})=\omega(\gamma_{24})=-1$ and $\omega(\gamma)=0$ 
if $\gamma$ is a $2$-cycle. 
Then it is easy to see that $\omega$ is weakly balanced 
on any pair of adjacent edges of $H$. 
We can see that 
$g_{1}$ contains exactly one 
non-trivial $2$-component constituent link 
$M=g_{1}(\gamma_{14}\cup \gamma'_{14})$ 
with ${\rm lk}(M)=0$ and $a_{3}(M)=2$. 
Thus by Theorem \ref{inv1} (2) we have that 
$\beta_{\omega,\omega}(g_{1})$ is an integer-valued vertex-homotopy 
invariant of $g_{1}$ and $\beta_{\omega,\omega}(g_{1})=2$. This implies 
that $g_{1}$ is non-splittable up to vertex-homotopy.  

\begin{figure}[htbp]
      \begin{center}
\scalebox{0.35}{\includegraphics*{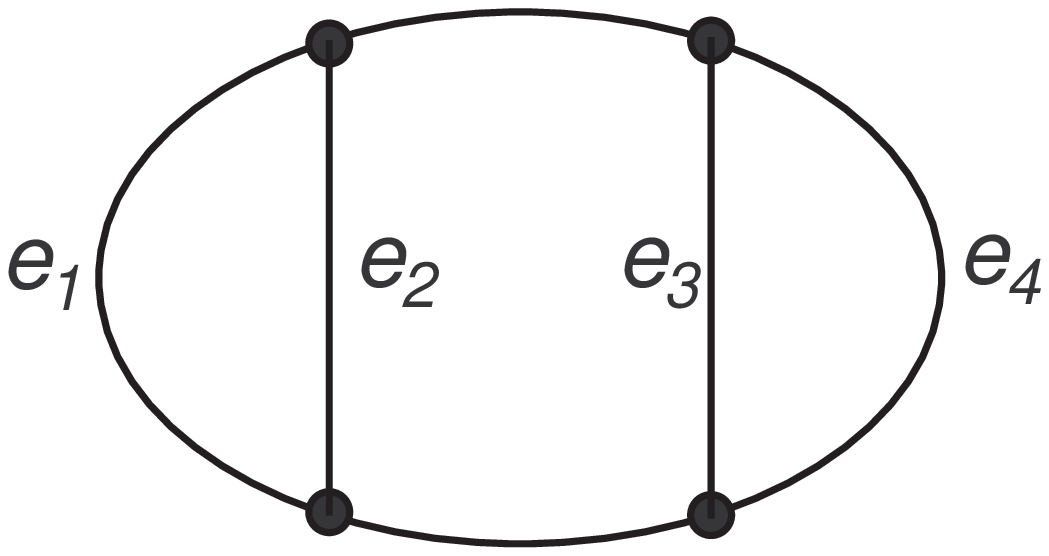}}
      \end{center}
   \caption{}
  \label{VH_ex_graph}
\end{figure} 
\begin{figure}[htbp]
      \begin{center}
\scalebox{0.33}{\includegraphics*{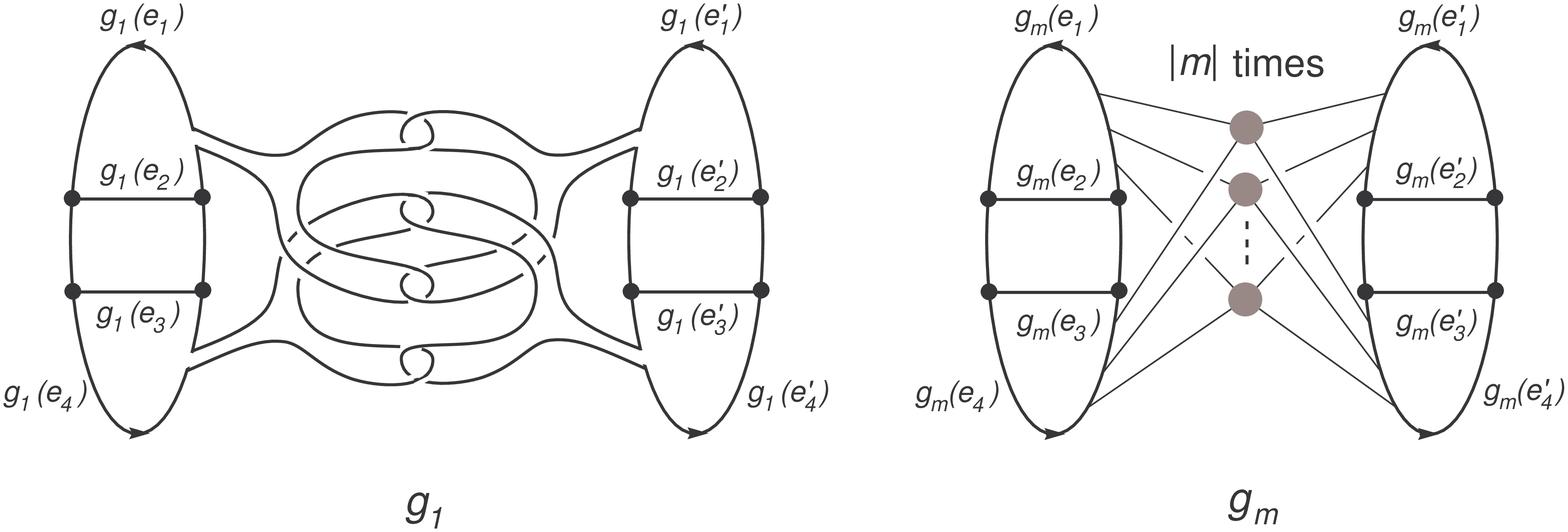}}
      \end{center}
   \caption{}
  \label{VH_ex}
\end{figure} 

Moreover, let $g_{m}$ be a spatial embedding of $G$ 
as illustrated in Figure \ref{VH_ex}, which can be constructed 
in the same way as in Example \ref{integer_inv_ex}. 
Then we can see that $\beta_{\omega,\omega}(g_{m})$ is an integer-valued 
vertex-homotopy invariant of $g_{m}$ and $\beta_{\omega,\omega}(g_{m})=2m$. 
This implies that $g_{m}$ is non-splittable up to vertex-homotopy for 
any integer $m\neq 0$ and $g_{i}$ and $g_{j}$ are not 
vertex-homotopic for any $i\neq j$.  
}
\end{Example}
\section{Modulo two invariants} 
\begin{Proposition}\label{modulo2_inv_weight}
Let $G$ be a planar graph which is not homeomorphic to $S^{1}$ 
and $p:G\to S^{2}$ a cellular embedding. 
Let $\omega_{p}:\Gamma(G)\to {\mathbb Z}_{2}$ be a weight on $\Gamma(G)$ over 
${\mathbb Z}_{2}$ defined by 
\begin{eqnarray*}
\omega_{p}(\gamma)= \left\{
\begin{array}{@{\,}ll}
1 & \mbox{$(\gamma\in \Gamma_{p}(G))$,} \\
0 &   \mbox{$(\gamma\in \Gamma(G)\setminus \Gamma_{p}(G))$.}
\end{array}
\right. 
\end{eqnarray*}
Then $\omega_{p}$ is totally balanced. 
\end{Proposition}

\begin{proof}
It holds that 
\begin{eqnarray*}
\sum_{\gamma\in \Gamma(G)}\omega_{p}(\gamma)[\gamma]
=\sum_{\gamma\in \Gamma_{p}(G)}[\gamma]
=2\left[\sum_{e\in E(G)}e\right]=0
\end{eqnarray*}
in $H_{1}(G;{\mathbb Z}_{2})$, where $E(G)$ denotes the set of all edges of $G$. 
Thus we have the result. 
\end{proof}

Thus by 
Proposition \ref{modulo2_inv_weight} and Theorem \ref{mod2inv} (1), 
we can obtain an edge-homotopy invariant as follows.

\begin{figure}[htbp]
      \begin{center}
\scalebox{0.35}{\includegraphics*{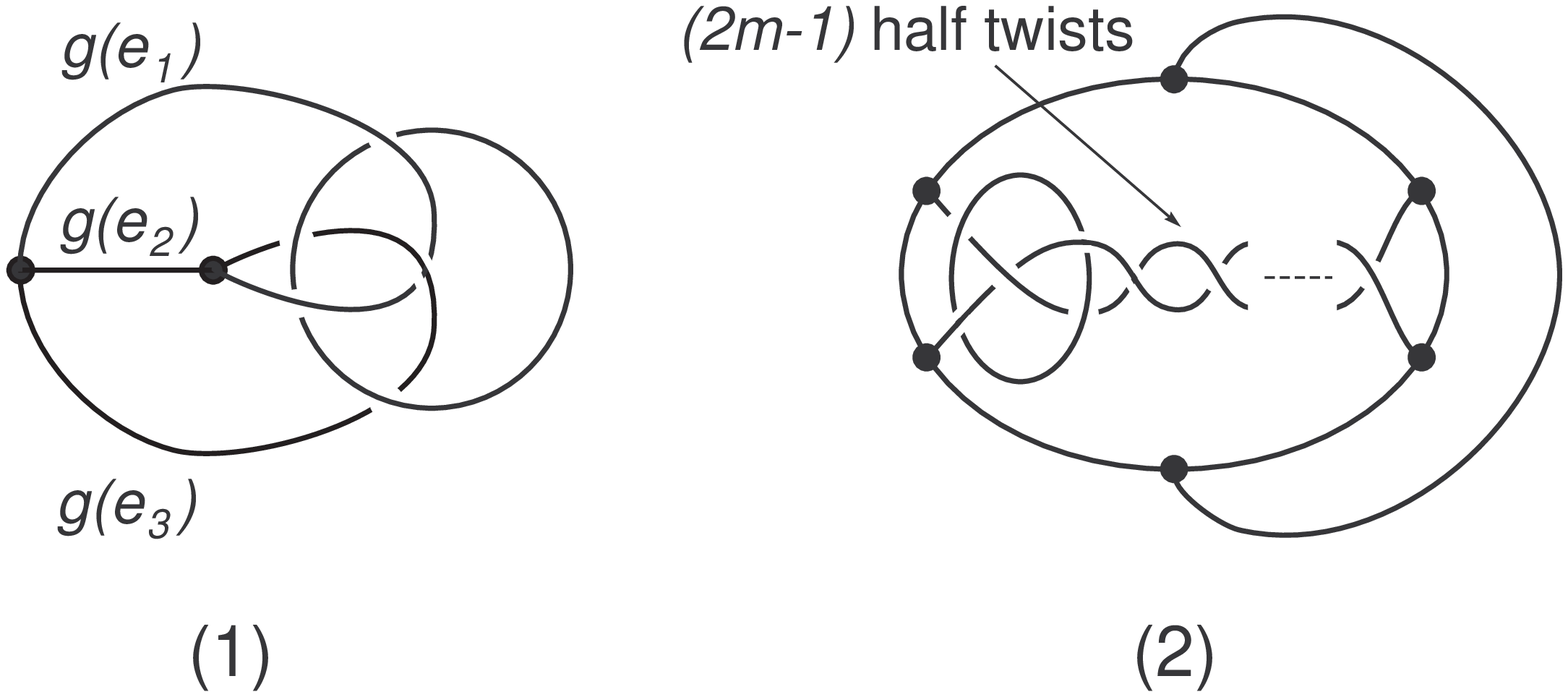}}
      \end{center}
   \caption{}
  \label{theta3_K33_circle}
\end{figure}

\begin{Theorem}\label{modulo2_inv}
Let $G=G_{1}\cup G_{2}$ be a disjoint union of two connected graphs such 
that $G_{1}$ is planar, not homeomorphic to $S^{1}$ and 
admits a cellular embedding $p_{1}:G_{1}\to S^{2}$.  
Let $\omega_{p_{1}}$ be a weight on $\Gamma(G_{1})$ over ${\mathbb Z}_{2}$ 
as in Proposition \ref{modulo2_inv_weight}, $\omega_{2}$ a weight 
on $\Gamma(G_{2})$ over ${\mathbb Z}_{2}$ and  
$f$ a spatial embedding of $G$ such that 
\begin{eqnarray*}
\omega_{p_{1}}(\gamma)\omega_{2}(\gamma')
{\rm lk}(f(\gamma),f(\gamma'))=0
\end{eqnarray*}
in ${\mathbb Z}$ for any $\gamma\in \Gamma(G_{1})$ and 
$\gamma'\in \Gamma(G_{2})$. 
Then $\beta_{\omega_{p_{1}},\omega_{2}}(f)$ is an edge-homotopy invariant 
of $f$. \hfill $\square$
\end{Theorem}
\begin{Example}\label{mod2inv_ex}
{\rm Let $G$ be a disjoint union of $\Theta_{3}$ and a circle $\gamma$. 
Let $\omega_{p}$ be a weight on $\Gamma(\Theta_{3})$ over ${\mathbb Z}_{2}$ 
as in Proposition \ref{modulo2_inv_weight} with respect to a cellular 
embedding $p:\Theta_{3}\to S^{2}$ as in Example \ref{integer_inv_ex},  
and $\omega$ a weight on 
$\Gamma(\gamma)$ over ${\mathbb Z}_{2}$ defined by $\omega(\gamma)=1$. 
Let $g$ be a spatial embedding of $G$ as illustrated in 
Figure \ref{theta3_K33_circle} (1). We can see that 
$g$ contains exactly one 
non-trivial $2$-component link 
$L=g(\gamma_{13})\cup g(\gamma)$ which is the Whitehead link, 
so ${\rm lk}(L)=0$ and $a_{3}(L)=1$.
Thus by Theorem \ref{modulo2_inv} we have that 
$\beta_{\omega_{p},\omega}(g)$ is an edge-homotopy 
invariant of $g$ and $\beta_{\omega_{p},\omega}(g)=1$.  
Namely $g$ is non-splittable up to edge-homotopy despite the fact that any of 
the constituent links of $g$ is link-homotopically trivial. 
}
\end{Example}
\begin{Example}\label{mod2inv_ex2}
{\rm 
Let $G$ be a disjoint union of the complete bipartite graph on $3+3$ vertices 
$K_{3,3}$ and a circle $\gamma$. 
Let $\omega_{3,3}$ be a weight on $K_{3,3}$ over ${\mathbb Z}_{2}$ 
defined by $\omega_{3,3}(\gamma')=1$ if $\gamma'$ is a $4$-cycle 
and $0$ if $\gamma'$ is a $6$-cycle. Let 
$\omega$ be a weight on 
$\Gamma(\gamma)$ over ${\mathbb Z}_{2}$ defined by $\omega(\gamma)=1$. 
Then it is not hard to see that $\omega_{3,3}$ is totally balanced and 
weakly balanced on any pair of adjacent edges of $K_{3,3}$. 
For a positive integer $m$, let $g_{m}$ be a spatial embedding of $G$ 
as illustrated in Figure \ref{theta3_K33_circle} (2). 
Note that $g_{i}(K_{3,3})$ and $g_{j}(K_{3,3})$ are not 
vertex-homotopic for any $i\neq j$ \cite{motohashi-taniyama96}; 
namely $g_{i}$ and $g_{j}$ are not vertex-homotopic for any $i\neq j$.  
Since all of the $2$-component constituent links of 
$g_{m}$ are algebraically split, by Theorem \ref{mod2inv} (2) 
we have that $\beta_{\omega_{3,3},\omega}(g)$ is a vertex-homotopy 
invariant of $g_{m}$. Moreover we can see that there exists exactly 
one $4$-cycle $\gamma'$ of $K_{3,3}$ so that 
$L=g_{m}(\gamma\cup\gamma')$ is non-trivial. Since $L$ is the Whitehead link,  
we have that $\beta_{\omega_{3,3},\omega}(g_{m})=1$. Therefore $g_{m}$ is 
non-splittable up to vertex-homotopy despite the fact that any of 
the constituent links of $g$ is link-homotopically trivial. 
}
\end{Example}
\begin{Remark}\label{mod2inv_rem}
{\rm 
The ${\mathbb Z}_{2}$-valued invariant in Theorem \ref{mod2inv} cannot always 
be extended to an integer-valued one. For example, 
\begin{enumerate}
\item Let us consider the graph $G$ and the invariant 
$\beta_{\omega_{p},\omega}$ as in Example \ref{mod2inv_ex}. 
Let $f$ be a spatial embedding of $G$ as illustrated in 
Figure \ref{not_int}. We can see that $f$ is edge-homotopic to the trivial 
spatial embedding $h$ of $G$. But by a calculation we have that 
$\sum_{1\le i<j\le 3}a_{3}(f(\gamma_{ij}),f(\gamma))=-2$. 

\begin{figure}[htbp]
      \begin{center}
\scalebox{0.35}{\includegraphics*{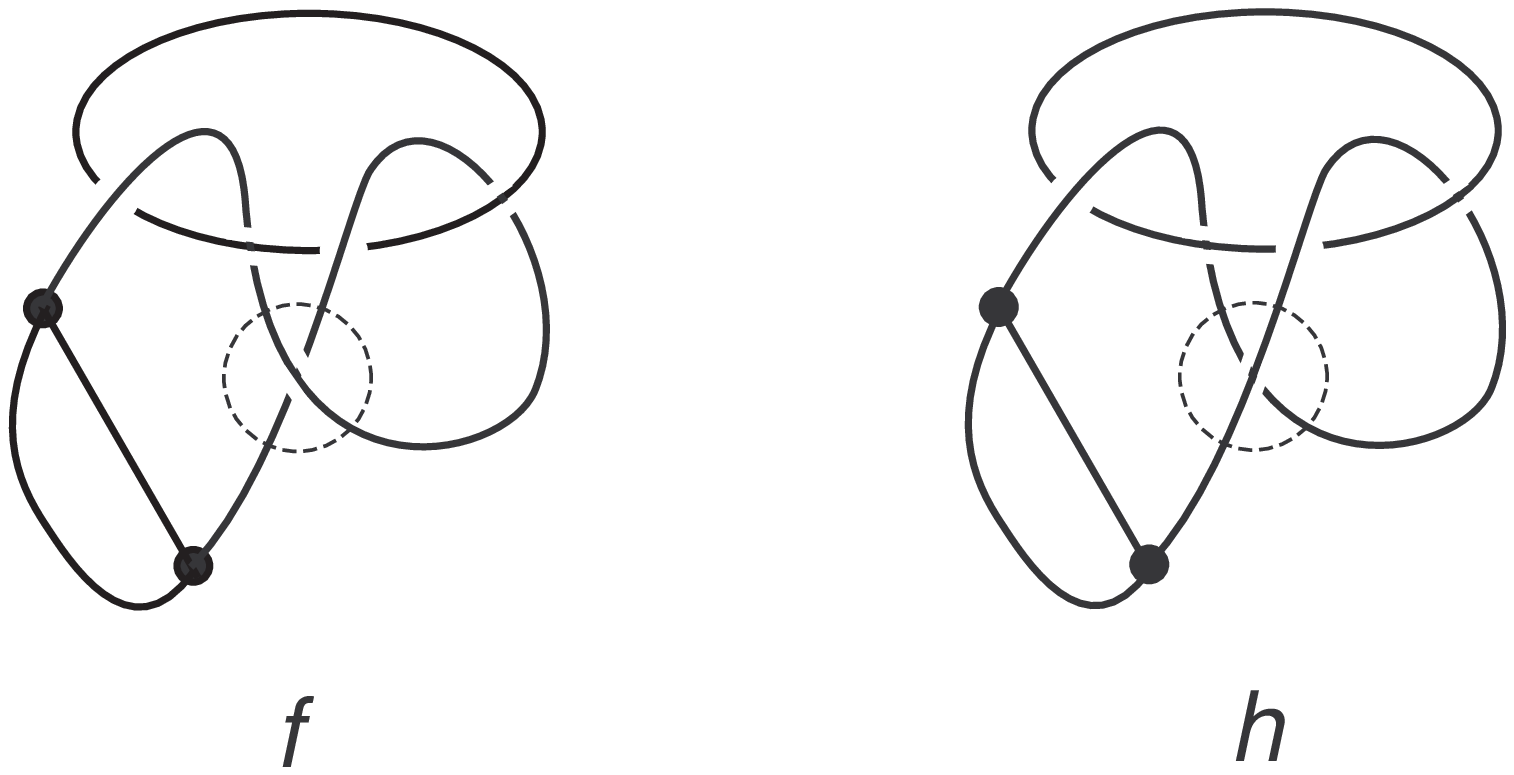}}
      \end{center}
   \caption{}
  \label{not_int}
\end{figure} 

\item Let $G$ be a disjoint union of $\Theta_{4}$ and a circle $\gamma$. 
Let $\omega_{p}$ be a checkerboard weight on $\Gamma(\Theta_{4})$ 
over ${\mathbb Z}$ 
as in Example \ref{integer_inv_ex}. 
Note that the modulo two reduction of a checkerboard weight is totally 
balanced. So by Theorem \ref{mod2inv} (1), the modulo two reduction of 
$\sum_{\gamma_{ij}\in\Gamma(\Theta_{4})}\omega_{p}(\gamma_{ij})
a_{3}(f(\gamma_{ij}\cup \gamma))$ is an edge-homotopy invariant of a 
spatial embedding $f$ of $G$. Moreover, we can see that 
the integer-value 
$\sum_{\gamma_{ij}\in\Gamma(\Theta_{4})}\omega_{p}(\gamma_{ij})
a_{3}(f(\gamma_{ij}\cup \gamma))$ is invariant under the self crossing 
change on $f(\Theta_{4})$ in the same way as in the proof of 
Theorem \ref{inv1} (1). But this value may change under a self crossing 
change on $f(\gamma)$. For example, let $f$ and $g$ be two spatial embeddings 
of $G$ as illustrated in 
Figure \ref{not_inv2}. We can see that $f$ is edge-homotopic to $g$. 
But by a calculation we have that 
\begin{eqnarray*}
\sum_{\gamma_{ij}\in \Gamma(\Theta_{4})}\omega_{p}(\gamma_{ij})
a_{3}(f(\gamma_{ij}),f(\gamma))&=&-1, \\ 
\sum_{\gamma_{ij}\in \Gamma(\Theta_{4})}\omega_{p}(\gamma_{ij})
a_{3}(g(\gamma_{ij}),g(\gamma))&=&1. 
\end{eqnarray*}
\end{enumerate}

\begin{figure}[htbp]
      \begin{center}
\scalebox{0.325}{\includegraphics*{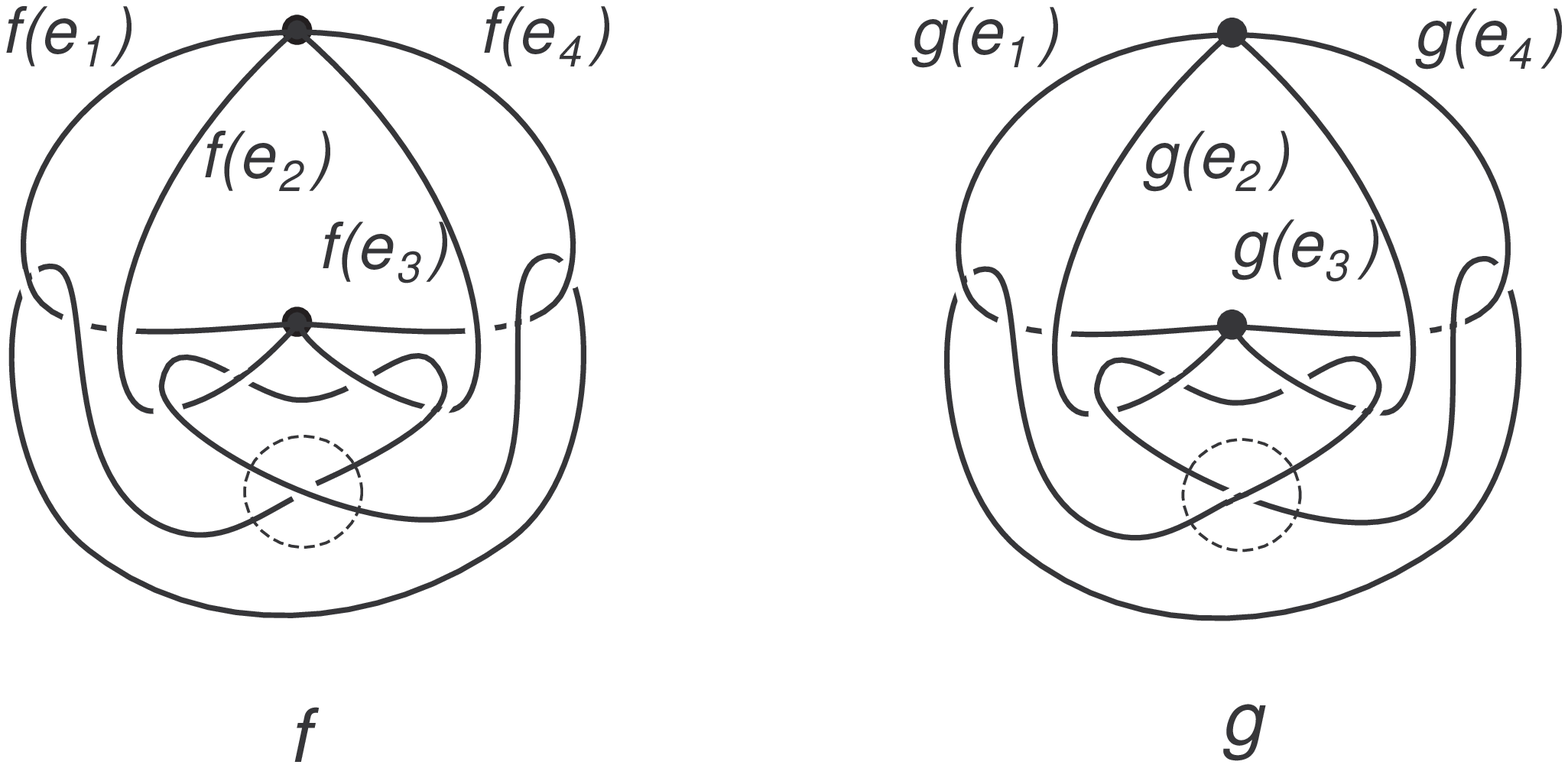}}
      \end{center}
   \caption{}
  \label{not_inv2}
\end{figure} 
}
\end{Remark}
\section{Applying the boundary of a planar surface} 

Let $X$ be a disjoint union of a graph $G$ and a planar surface $F$ 
with boundary. 
Let $\omega$ be a weight on $\Gamma(G)$ over ${\mathbb Z}_{2}$ and 
$\varphi$ an embedding of $X$ into $S^{3}$ 
such that 
\begin{eqnarray*}
\omega(\gamma)
{\rm lk}(\varphi(\gamma),\varphi(\gamma'))=0
\end{eqnarray*}
in ${\mathbb Z}$ for any $\gamma\in \Gamma(G)$ and 
$\gamma'\in \Gamma(\partial F)$. 
Then we define $\beta_{\omega}(\varphi)\in {\mathbb Z}_{2}$ by 
\begin{eqnarray*}
\beta_{\omega}(\varphi)\equiv \sum_{\gamma\in\Gamma(G) \atop \gamma'
\in\Gamma(\partial F)}
\omega(\gamma)
a_{3}(\varphi(\gamma),\varphi(\gamma'))\pmod{2}. 
\end{eqnarray*}
Let $G$ be a disjoint union of a connected graph $G_{1}$ and 
a connected planar graph $G_{2}$. 
Let $f$ be a spatial embedding of $G$ and $p$ an embedding 
of $G_{2}$ into $S^{2}$. We denote the regular neighborhood of 
$p(G_{2})$ in $S^{2}$ by $F(G_{2};p)$, which is a planar surface having 
$p(G_{2})$ as a spine. Then the spatial embedding $f$ 
induces an embedding $\tilde{f}_{p}$ of the disjoint union 
$G_{1}\cup F(G_{2};p)$ 
into $S^{3}$, 
so that $\tilde{f}_{p}(G_{1})=f(G_{1})$ and 
$\tilde{f}_{p}(F(G_{2};p))$ has $f(G_{2})$ as a spine in the natural way. 
Note that such an induced embedding $\tilde{f}_{p}$ is not unique 
up to ambient isotopy. 
Let $\omega$ be a weight on $\Gamma(G_{1})$ over ${\mathbb Z}_{2}$ 
so that 
\begin{eqnarray*}
\omega(\gamma)
{\rm lk}(\tilde{f}_{p}(\gamma),\tilde{f}_{p}(\gamma'))=0
\end{eqnarray*}
in ${\mathbb Z}$ for any $\gamma\in \Gamma(G_{1})$ and 
$\gamma'\in \Gamma(\partial F(G_{2};p))$. Then we have the following. 

\begin{Theorem}\label{modulo2_inv_more}
If $f$ is split up to edge-homotopy, then 
$\beta_{\omega}(\tilde{f}_{p})=0$ for any induced embedding 
$\tilde{f}_{p}$ of $G_{1}\cup F(G_{2};p)$. 
\end{Theorem}

\begin{figure}[htbp]
      \begin{center}
      \scalebox{0.35}{\includegraphics*{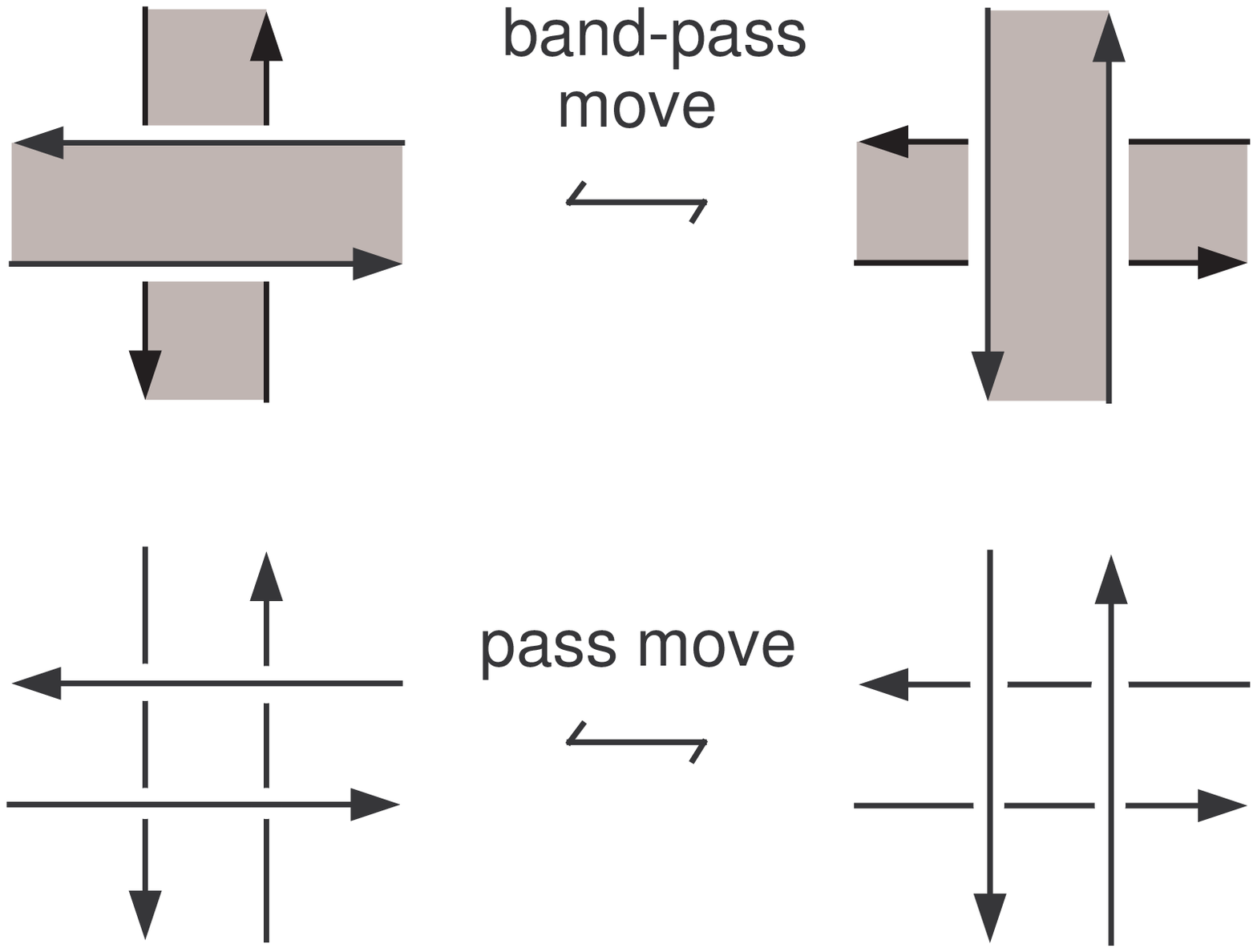}}
\end{center}
   \caption{}
  \label{band-pass}
\end{figure}

\begin{figure}[htbp]
      \begin{center}
\scalebox{0.325}{\includegraphics*{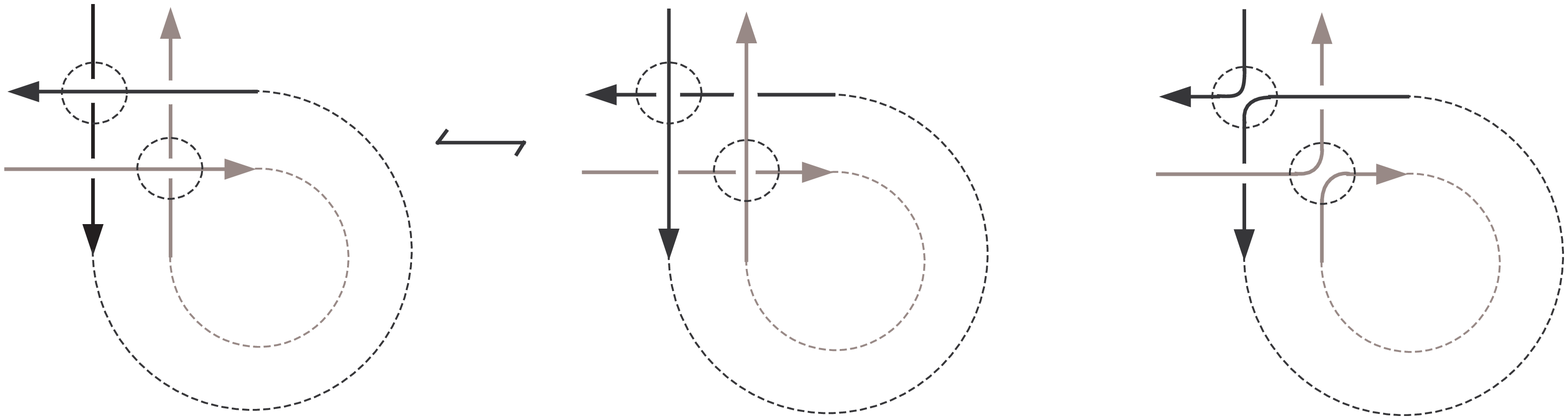}}
      \end{center}
   \caption{}
  \label{pass-smoothing}
\end{figure}

\begin{proof}
By the assumption we have that $f$ is transformed into 
a split spatial embedding $u$ of $G$ by self crossing changes and 
ambient isotopies. Then each of the self crossing changes induces a 
self crossing change on $\tilde{f}_{p}(G_{1})$ or a 
{\it band-pass move} \cite{kauffman83} (see Figure \ref{band-pass}) on 
$\tilde{f}_{p}(F(G_{2};p))$. Namely $\tilde{f}_{p}$ can be transformed into 
an induced embedding $\tilde{u}_{p}$ of $G_{1}\cup F(G_{2};p)$ by 
such moves and ambient isotopies. 
Let $\tilde{g}_{p}$ be an embedding of 
$G_{1}\cup F(G_{2};p)$ into $S^{3}$ obtained from $\tilde{f}_{p}$ 
by a single self crossing change on $\tilde{f}_{p}(G_{1})$ or a single 
band-pass move on $\tilde{f}_{p}(F(G_{2};p))$. 
Then it still holds that 
\begin{eqnarray*}
\omega(\gamma){\rm lk}(\tilde{g}_{p}(\gamma),\tilde{g}_{p}(\gamma'))=0
\end{eqnarray*}
in ${\mathbb Z}$ for any $\gamma\in \Gamma(G_{1})$ and 
$\gamma'\in \Gamma(\partial F(G_{2};p))$.  
\begin{flushleft}
{\bf Claim.} $\beta_{\omega}(\tilde{f}_{p})=\beta_{\omega}(\tilde{g}_{p})$. 
\end{flushleft}

Assume that $\tilde{g}_{p}$ is obtained from $\tilde{f}_{p}$ 
by a single self crossing change on $\tilde{f}_{p}(G_{1})$. 
Since it holds that 
\begin{eqnarray*}
\sum_{\gamma'\in \Gamma(\partial F(G_{2};p))}[\gamma']=0
\end{eqnarray*}
in $H_{1}(F(G_{2};p);{\mathbb Z}_{2})$, 
we can see that 
$\beta_{\omega}(\tilde{f}_{p})=\beta_{\omega}(\tilde{g}_{p})$ 
in a similar way as the proof of 
Theorem \ref{mod2inv} (1). Next we assume that $\tilde{g}_{p}$ 
is obtained from $\tilde{f}_{p}$ 
by a single band-pass move on $\tilde{f}_{p}(F(G_{2};p))$. 
Then $\tilde{g}_{p}|_{G_{1}\cup \partial F(G_{2};p)}$ is obtained 
from $\tilde{f}_{p}|_{G_{1}\cup \partial F(G_{2};p)}$ by a single 
{\it pass move} \cite{kauffman83} (see Figure \ref{band-pass}) 
on $\tilde{f}_{p}(\partial F(G_{2};p))$. We 
divide our situation into the following two cases. 

\vspace{0.2cm}
{\bf Case 1.} Four strings in the pass move belong to 
$\tilde{f}_{p}(\gamma_{1}')$ and $\tilde{f}_{p}(\gamma_{2}')$ for 
exactly two cycles $\gamma_{1}'$ and $\gamma_{2}'$ in 
$\Gamma(\partial F(G_{2};p))$. 

\vspace{0.2cm}
This pass move causes a single self crossing change on 
$\tilde{f}_{p}(\gamma_{1}')$ and a single self crossing change on 
$\tilde{f}_{p}(\gamma_{2}')$. Then the separated components that result from 
smoothing each of the self crossings are orientation-reversing parallel 
knots; see Figure \ref{pass-smoothing}. 
So the difference between $\beta_{\omega}(\tilde{f}_{p})$ and 
$\beta_{\omega}(\tilde{g}_{p})$ is  
cancelled out in a similar way as in the proof of Theorem \ref{inv1} (1). 
Thus we have that 
$\beta_{\omega}(\tilde{f}_{p})=\beta_{\omega}(\tilde{g}_{p})$. 

\vspace{0.2cm}
{\bf Case 2.} Four strings in the pass move belong to 
$\tilde{f}_{p}(\gamma')$ for a cycle $\gamma'$ in 
$\Gamma(\partial F(G_{2};p))$. 

\vspace{0.2cm}
It is known that a pass move on the same component of a proper link 
$L=J_{1}\cup J_{2}\cup \cdots\cup J_{n}$ preserves 
$\overline{{\rm Arf}}(L)\equiv {\rm Arf}(L)-\sum_{i=1}^{n}{\rm Arf}(J_{i})\in 
{\mathbb Z}_{2}$ (cf. \cite{shibuya-yasuhara03}).\footnote{The value of 
$\overline{{\rm Arf}}(L)$ is called the {\it reduced Arf invariant} 
of $L$ \cite{shibuya89}.} 
Especially, if $n=2$ then $a_{3}(L)\equiv \overline{{\rm Arf}}(L)\pmod{2}$ \cite[Lemma 3.5 (ii)]{nikkuni05}. 
Therefore in this case the pass move preserves 
$\omega(\gamma)a_{3}(\tilde{f}_{p}(\gamma),\tilde{f}_{p}(\gamma'))$ 
for any cycle $\gamma\in \Gamma(G_{1})$. This implies that 
$\beta_{\omega}(\tilde{f}_{p})=\beta_{\omega}(\tilde{g}_{p})$. 

\vspace{0.2cm}
Now by the argument above, 
we have that 
$\beta_{\omega}(\tilde{f}_{p})=\beta_{\omega}(\tilde{u}_{p})$. 
Then, each $2$-component link 
$\tilde{u}_{p}(\gamma\cup\gamma')$ is split 
for any $\gamma\in\Gamma(G_{1})$ and 
$\gamma'\in\Gamma(\partial F(G_{2};p))$ 
because $u$ is split. 
Therefore we have that 
$\beta_{\omega}(\tilde{f}_{p})=\beta_{\omega}(\tilde{u}_{p})=0$. 
This completes the proof. 
\end{proof}

\begin{Example}\label{hand_bouquet_ex}
{\rm Let $G$ be a disjoint union of a circle $\gamma$ and 
the {\it handcuff graph} (resp. {\it $2$-bouquet}) $G_{2}$. 
Let $\omega$ be a weight on $\Gamma(\gamma)$ over ${\mathbb Z}_{2}$ 
defined by $\omega(\gamma)=1$. 
We fix an embedding $p:G_{2}\to S^{2}$ 
and take a regular neighborhood $F(G_{2};p)$ as illustrated in 
Figure \ref{nbd_plane} (1) (resp. (2)). 

\begin{figure}[htbp]
      \begin{center}
\scalebox{0.35}{\includegraphics*{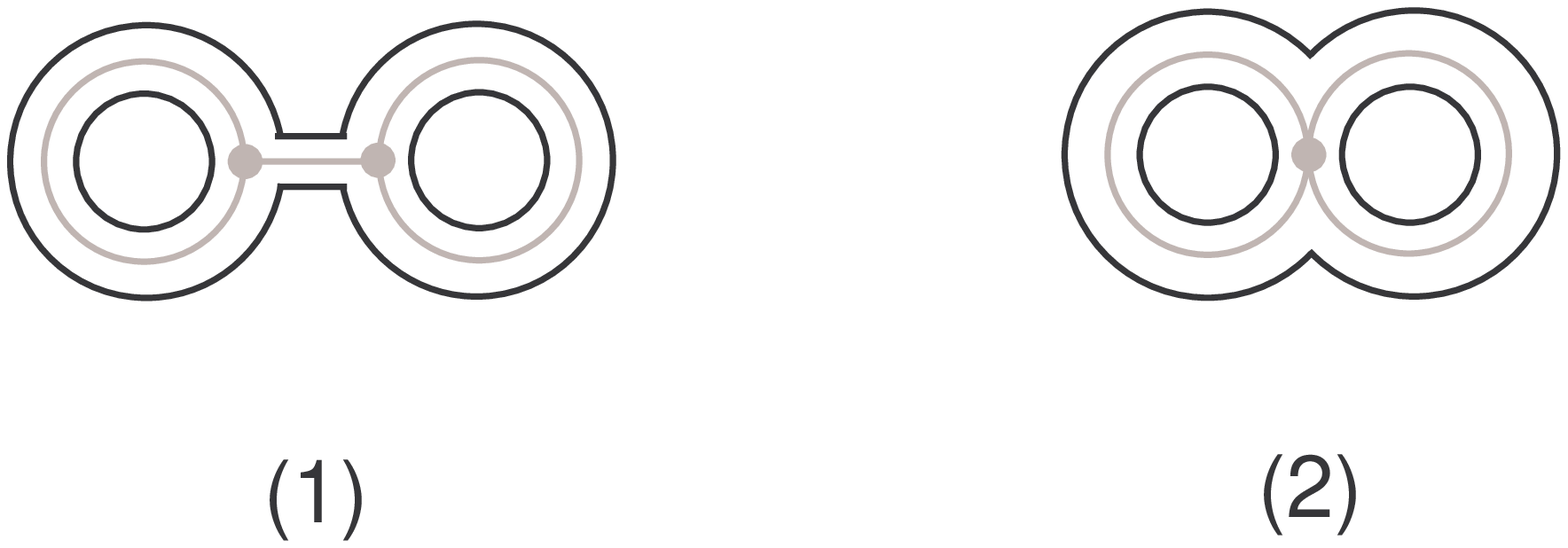}}
\end{center}
   \caption{}
  \label{nbd_plane}
\end{figure} 

Let $f$ be a spatial embedding of $G$ as illustrated in 
Figure \ref{handcuff_bouquet} (1) (resp. (2)). Let us take an induced 
embedding $\tilde{f}_{p}:\gamma\cup F(G_{2};p)\to S^{3}$ as illustrated 
in Figure \ref{handcuff_bouquet2} (1) (resp. (2)). Note that  
${\rm lk}(\tilde{f}_{p}(\gamma),\tilde{f}_{p}(\gamma'))=0$ for any 
$\gamma'\in \Gamma(\partial F(G_{2};p))$. 
Then it can be calculated that $\beta_{\omega}(\tilde{f}_{p})=1$. 
Thus by Theorem \ref{modulo2_inv_more} we have that $f$ is non-splittable 
up to edge-homotopy. 

\begin{figure}[htbp]
      \begin{center}
\scalebox{0.35}{\includegraphics*{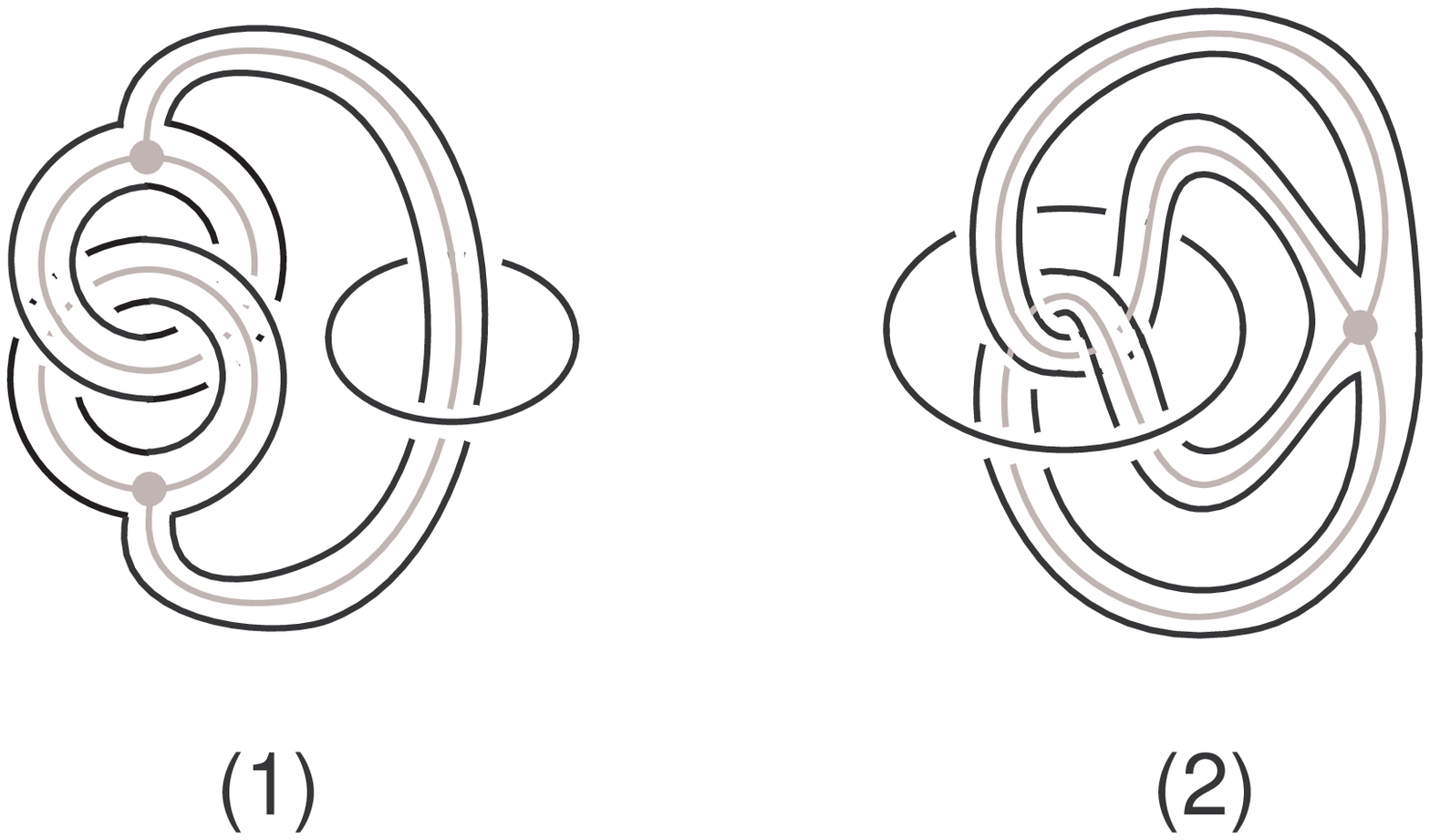}}
      \end{center}
   \caption{}
  \label{handcuff_bouquet2}
\end{figure} 
}
\end{Example}

\section*{Acknowledgment}

The authors are very grateful to Professor Hitoshi Murakami for 
his hospitality at the Tokyo Institute of Technology where this work 
was conducted. 

{\normalsize
}

\end{document}